\documentclass[reprint,nofootinbib,secnumarabic,amssymb, nobibnotes, aps, prx]{revtex4-2}

\usepackage[utf8]{inputenc}
\usepackage{amsmath}
\usepackage{amsfonts}
\usepackage{pmat}
\usepackage{amssymb}
\usepackage{amsthm}
\usepackage{bm}
\usepackage{xcolor}
\usepackage{fullpage}
\usepackage{enumitem}
\usepackage{comment}
\usepackage{bbm}
\usepackage{tikz}
\usepackage{graphicx}
\usepackage{circuitikz}
\usepackage{pgfplots}
\usepgfplotslibrary{fillbetween}
\pgfplotsset{compat=1.5}
\newcommand{\norm}[1]{\left\lVert#1\right\rVert}

    \usetikzlibrary{shapes,arrows}
    \usetikzlibrary{arrows,calc,positioning}
\tikzset{%
        block/.style = {draw, rectangle,
            minimum height=1cm,
            minimum width=2cm},
        input/.style = {coordinate,node distance=1.5cm},
        output/.style = {coordinate,node distance=1cm},
        arrow/.style={draw, -latex,node distance=1cm},
        pinstyle/.style = {pin edge={latex-, black,node distance=2cm}},
	system/.style = {rectangle, draw=black, line width=0.5pt, fill=rugblue!30,
		inner sep=0pt, minimum width=14mm, minimum height=8mm,
		text height=1.8ex,text depth=.25ex},
	sum/.style = {circle, draw=black, line width=0.5pt, fill=white, inner sep=0pt, minimum size=0.4cm},
	signal/.style = {color=black, line width=0.7pt, >=latex},
}

\newcommand{\ones}{\mathbbm{1}}

\newcommand{\R}{\mathbb{R}}
\newcommand{\N}{\mathbb{N}}
\newcommand{\calC}{\mathcal{C}}
\newcommand{\calV}{\mathcal{V}}
\newcommand{\calG}{\mathcal{G}}
\newcommand{\calE}{\mathcal{E}}

\newcommand{\calX}{\mathcal{X}}
\newcommand{\calY}{\mathcal{Y}}
\newcommand{\calF}{\mathcal{F}}
\newcommand{\calA}{\mathcal{A}}

\DeclareMathOperator{\im}{im}

\newtheorem{theorem}{Theorem}
\newtheorem{lemma}[theorem]{Lemma}

\newtheorem{proposition}[theorem]{Proposition}
\newtheorem{definition}[theorem]{Definition}
\newtheorem{example}[theorem]{Example}

\newtheorem{remark}[theorem]{Remark}
\newtheorem{assumption}[]{Assumption}

\DeclareMathOperator*{\argmin}{arg\,min}
\DeclareMathOperator*{\diag}{diag}
\DeclareMathOperator{\smax}{smax}
\DeclareMathOperator{\smin}{smin}
\DeclareMathOperator*{\inter}{int}

\begin{document}

\title{Target-Clamped Contrastive Learning in Nonlinear Resistive Networks \\ with Monotone Elements}

\author{Henk J. van Waarde$^{1}$}
\email{Corresponding author: h.j.van.waarde@rug.nl}

\author{Benjamin Scellier}

\author{Encho T. Nguyen$^{1}$}

\affiliation{$^{1}$Bernoulli Institute for Mathematics, Computer Science and AI, University of Groningen, Nijenborgh 9, 9747AG, Groningen, The Netherlands}
\date{\today}


\begin{abstract}
In the pursuit of energy-efficient hardware for machine learning, nonlinear resistive networks trained via local learning rules are emerging as promising alternatives to digital neural networks. Learning in these physical systems is often achieved with Equilibrium Propagation (EP), which compares a ``free" state to a ``nudged" target-driven state to perform gradient descent on an objective defined by an \emph{explicit} cost function. In contrast, Target-Clamped Contrastive Learning (TCCL) ``clamps" output nodes directly to exact target values, performing gradient descent on an \emph{implicit} contrastive function---an objective whose properties have long remained unknown. Focusing on nonlinear networks composed of resistors with monotonically increasing current-voltage characteristics (such as Ohmic resistors and diodes), our contribution is twofold. First, we establish rigorous convergence guarantees for TCCL. We prove that the contrastive function is convex for nonlinear resistive networks without hidden nodes, and demonstrate that this property does not hold in general for the objective function optimized by EP. For these networks, TCCL converges to a global minimum of the contrastive function under suitable assumptions on the learning rate. In the presence of hidden nodes, the contrastive function is not convex, but we show that TCCL still converges to a stationary point on a subsequence. We also reveal cases where such stationary points are global minimizers of the contrastive function. Second, we characterize the representational capacity of these circuits. We establish a novel function approximation theorem demonstrating that any continuous, monotonically non-decreasing, and shift-equivariant function can be approximated to arbitrary precision, proving that these properties are the sole limitations on the network's computational expressivity.                   
\end{abstract}

\maketitle

\section{Introduction}

The escalating energy demands of modern machine learning have driven a resurgence of interest in analog computing paradigms and the development of training methods dedicated to such platforms \citep{momeni2025training}. By replacing the digital von Neumann architecture with analog compute-in-memory substrates, these systems promise substantial improvements in energy efficiency for both inference and training. Among these alternative platforms, resistive networks have attracted special interest.

A resistive network is an electrical circuit comprised of (linear or nonlinear) resistors, with nodes typically partitioned into input, hidden and output sets. In such a network where computation is driven by the circuit's physics, the input node potentials are clamped to specific values, while the output node potentials serve as the network's prediction. Conductances of variable resistors serve as the trainable parameters, and the objective is to find parameter configurations that implement a desired input-output mapping. It is well-known that the hidden and output potentials of a linear resistive network minimize the total dissipated power; in nonlinear resistive networks---those containing nonlinear one-port resistive elements---the hidden and output potentials extremize a more general notion known as the co-content \citep{Millar1951}, also referred to as “pseudo-power” \citep{kendall2020training}. This property allows resistive networks to function as energy-based models where the co-content acts as the energy function. For such systems, energy-based learning rules have been developed which adjust conductances based solely on local voltage measurements, allowing decentralized training  \citep{kendall2020training,dillavou2022demonstration,anisetti2024frequency,huijzer2026convergence}. On the hardware side, small-scale prototypes have successfully demonstrated self-learning capabilities on toy tasks \citep{dillavou2022demonstration,dillavou2024machine}. On the simulation side, SPICE-class simulations of modest-size nonlinear resistive networks have been performed on MNIST \citep{kendall2020training,watfa2023energy}, while recent methods exploit the co-content minimization formulation to achieve speedups of several orders of magnitude \citep{pmlr-v235-scellier24a,osana2026fast}, paving the way for large-scale simulations.

To rival artificial neural networks, resistive networks must demonstrate representational capacity, and learning capabilities. Learning in resistive networks was initially formalized by \cite{kendall2020training} through Equilibrium Propagation (EP), an energy-based learning framework that enables gradient descent on an arbitrary cost function \citep{scellier2017equilibrium}. EP is part of a broader class of energy-based learning methods \citep{Scellier2023} that optimize parameters by comparing the network in two different states: the inference state (called \emph{free} state), and a perturbed state where outputs are clamped or nudged towards a target. Other energy-based learning rules differ from EP in the way the perturbed state is defined. In particular, the Contrastive Learning (CL) algorithm of \cite{huijzer2026convergence} (adapted from \cite{Movellan1991,Scellier2023}) clamps output nodes to exact target values in the perturbed state, effectively performing gradient descent on a nonnegative \emph{contrastive function} whose minimum of zero is attained precisely when the network outputs match the target output. In \cite{huijzer2026convergence}, it was shown that the contrastive function is convex for linear networks without hidden nodes, establishing the convergence for a range of step sizes. These results were also recently extended to linear resistive-capacitive networks in \cite{Heidema2026}. On the representational capacity front, \cite{Scellier2025} showed that nonlinear resistive networks can approximate arbitrary continuous functions, provided they incorporate two-port elements like voltage-controlled voltage sources (VCVSs).

In this paper, we establish mathematical foundations for learning and representation in nonlinear resistive networks constructed \textit{exclusively} from one-port elements. This deliberately precludes the VCVSs utilized by \cite{Scellier2025}. We assume all circuit components possess strictly monotonically increasing current-voltage characteristics, which encompasses Ohmic resistors and diodes, and yields a strictly convex network co-content. Within this framework, we answer two fundamental questions: \emph{under what conditions does learning reliably converge, and what is the network's representational capacity?}

While \cite{huijzer2026convergence} proved the convergence of the CL algorithm with clamped outputs, their analysis was limited to linear networks without hidden nodes. Our first main contribution extends their framework to nonlinear networks. We prove that the contrastive function is convex for nonlinear networks without hidden nodes, guaranteeing convergence of the parameters to a global minimum. This is in contrast with the EP objective function, which is in general non-convex, even without hidden nodes. For the more expressive and general case that includes hidden nodes, the contrastive function is no longer convex, but in this case, we prove that the learning dynamics converge on a subsequence to a stationary point. Interestingly, any of such stationary points in the interior of the constraint set is shown to be a global minimum of the contrastive function.

Our second main contribution analyzes the representational expressivity of these circuits. We demonstrate that the input-output mappings of such networks are monotonically non-decreasing and shift-equivariant (translation-equivariant). Remarkably, we establish a novel approximation theorem proving that these two physical constraints are the only limitations on expressivity. Any continuous function satisfying monotonicity and shift equivariance can be approximated to arbitrary precision by a nonlinear resistive network featuring strictly convex co-content functions. To prove this, we introduce a novel circuit construction leveraging scalable, smooth approximations of min-max functions.

Finally, we clarify a terminological ambiguity present in the recent literature on physical learning algorithms. The term “contrastive learning” is used in two distinct ways: as the specific learning rule involving clamped targets \citep{Scellier2023,huijzer2026convergence,Heidema2026}, and as a broad umbrella term for the entire class of energy-based learning methods \citep{dillavou2022demonstration,Stern2023annualreview,momeni2025training,du2026metamaterials}. This broader category includes other algorithms such as coupled learning \citep{Stern2021,dillavou2022demonstration,Altman2024}, which employs a weaker form of clamping that fixes outputs to a convex combination of the free outputs and targets. To avoid confusion, we introduce the term \emph{Target-Clamped Contrastive Learning} (TCCL) to explicitly refer to the specific learning rule that utilizes clamped targets. For brevity, we occasionally abbreviate this to Contrastive Learning (CL) when the context is unambiguous, in line with the terminology of \citep{Scellier2023,huijzer2026convergence,Heidema2026}.

We note a recent concurrent study by \cite{lipka2026thermodynamic} on the computational expressivity of thermodynamic networks; however, their framework permits elements whose current depends on terminal voltages independently rather than on their difference, meaning the network non-equilibrium steady state cannot generally be characterized as the minimum of a co-content function, and energy-based learning algorithms such as TCCL or EP cannot be directly employed.

Another recent work \cite{mcginnis2026coercivity} studied the convergence of energy-based learning algorithms for linear resistive circuits. This work is complementary to the current paper, in the sense that it deals with local convergence results for EP and coupled learning in linear circuits, whereas our study focuses on \emph{global} convergence analysis for TCCL in \emph{nonlinear} circuits.

The remainder of this paper is organized as follows. Section \ref{sec:modeling} introduces nonlinear resistive circuits and co-content functions. Section \ref{sec:problem} formalizes the learning problem. Section~\ref{sec:CL} discusses the contrastive learning framework. Section~\ref{sec:convex} analyzes the convexity of the contrastive function, and Section~\ref{sec:algorithm} establishes convergence guarantees for the learning dynamics. Section~\ref{sec:approximation} defines the fundamental physical constraints of the circuits and proves the approximation theorem. Finally, conclusions are provided in Section~\ref{sec:conclusions}. The reader is also directed to Appendix~\ref{app:notation} for relevant notation, and to Appendix~\ref{app:convexity} for preliminaries on convex, coercive and Lipschitz continuous functions.

\section{Modeling of nonlinear resistive circuits}
\label{sec:modeling}

We start by introducing a few graph-theoretic concepts used throughout the paper. For a general reference on graph theory, we refer the reader to \cite{Bollobas1998}. We consider a connected undirected graph 
\begin{equation*}
\mathcal{G} = (\mathcal{V},\mathcal{E}),
\end{equation*}
consisting of a set of nodes $\mathcal{V} = \{1,2,\dots,N\}$ and a set of branches (or edges) $\mathcal{E}$, containing $B$ unordered pairs of nodes. The \emph{incidence matrix} of $\mathcal{G}$ is defined as follows. First, we denote the edges in the graph by $e_1,e_2,\dots,e_B$, and we assign an arbitrary orientation to each edge. Then, the incidence matrix $D \in \mathbb{R}^{N \times B}$ is defined element-wise by
\begin{equation*}
D_{ij} = \begin{cases}
    -1 & \text{if edge } e_j \text{ leaves node } i, \\
    1 & \text{if edge } e_j \text{ enters node } i, \\
    0 & \text{otherwise}.
\end{cases}
\end{equation*}
By definition of the incidence matrix, we have that $\im \ones \subseteq \ker D^\top$, where $\ones$ denotes the vector of ones. Moreover, since $\mathcal{G}$ is connected, it follows from \cite[Thm. 8.3.1]{Godsil2001} that
\begin{equation}
\label{kernelD}
\ker D^\top = \im \ones.
\end{equation}

Next, we consider the vectors $p \in \mathbb{R}^N$ of nodal voltage potentials, and $j \in \mathbb{R}^N$ of external currents injected into each of the nodes. The vectors $v \in \mathbb{R}^B$ and $i \in \mathbb{R}^B$ denote the voltages and currents over and through the edges, respectively. Obviously, due to the laws of Kirchhoff, we have 
\begin{equation}
\label{KCL-KVL}
v = D^\top p \:\: \text{ and }\:\: j = Di.
\end{equation}
Every branch $k \in \{1,2,\dots,B\}$ contains a \emph{nonlinear resistor} modeled as
\begin{equation}
\label{i-v-characteristic}
i_k = \frac{dG_k}{dv_k}(v_k).
\end{equation}
Here $G_k : \mathbb{R} \to \mathbb{R}$ is the so-called \emph{co-content function} \cite{Millar1951} of the $k$-th resistor, which we assume to be differentiable, coercive and strictly convex. Note that the strict convexity of $G_k$ is equivalent to $\frac{dG_k}{dv_k}$ being strictly monotonically increasing. Thus, \eqref{i-v-characteristic} defines a strictly monotone relation between the current and voltage at branch $k$.

\begin{remark}
\label{rem:examplesresistors}
Some simple examples of co-content functions are as follows. The case $G_k(v_k) = \frac{1}{2} g_k v_k^2$, where $g_k > 0$ is a conductance, covers linear resistors described by $i_k = g_k v_k$. If $G_k(v_k) = e^{v_k} - v_k$, then $i_k = e^{v_k} - 1$ is an (idealized) model of a p-n diode, see e.g., \cite{Chua1987}. Both of these examples of $G_k$ are differentiable, coercive and strictly convex. An example of a co-content function that does not satisfy the strict convexity assumption is that of a tunnel diode \cite[Ch. 2]{Desoer1969}, whose current-voltage characteristic is not monotonically increasing.
\end{remark} 

Next, we define the ``total" co-content function $G: \mathbb{R}^B \to \mathbb{R}$ as 
\begin{equation*}
G(v) := \sum_{k = 1}^B G_k(v_k).
\end{equation*}
The function $G$ inherits properties of the functions $G_k$. In particular, $G$ is differentiable, coercive and strictly convex. We also observe that 
\begin{equation*}
\nabla G(v) = \begin{bmatrix}
\frac{dG_1}{dv_1}(v_1) \\
\frac{dG_2}{dv_2}(v_2) \\ 
\vdots \\ 
\frac{dG_B}{dv_B}(v_B)
\end{bmatrix}.
\end{equation*}
This implies that, by Kirchhoff's laws \eqref{KCL-KVL} and the model of the resistors \eqref{i-v-characteristic}, we can succinctly write the relation between $j$ and $p$ as
\begin{equation}
\label{staticnetwork}
j = D \nabla G(D^\top p).
\end{equation}

In the context of learning in electrical circuits, it will become useful to partition the set of nodes $\mathcal{V}$ into two or multiple sets (corresponding to sets of input, hidden, and output nodes). In what follows, we will derive some preliminary results to aid that discussion. 

Let $\mathcal{V} = \mathcal{V}_1 \cup \mathcal{V}_2$ where $\mathcal{V}_1$ and $\mathcal{V}_2$ are disjoint. After possible relabeling of the nodes, this yields the partitioned vectors and matrix:
\begin{equation*}
p = \begin{bmatrix}
p_1 \\ p_2
\end{bmatrix}, \quad j = \begin{bmatrix}
j_1 \\ j_2
\end{bmatrix}, \quad D = \begin{bmatrix}
D_1 \\ D_2
\end{bmatrix},
\end{equation*}
where $p_1, j_2 \in \mathbb{R}^{N_1}$, $p_2,j_2 \in \mathbb{R}^{N_2}$, and $D_1 \in \mathbb{R}^{N_1 \times B}$, $D_2 \in \mathbb{R}^{N_2 \times B}$. Here $N_1,N_2 \in \mathbb{N}$ are positive integers such that $N_1 + N_2 = N$.

We set the injected currents in the $\mathcal{V}_2$-nodes to zero, i.e., $j_2 = 0$. By \eqref{staticnetwork}, this results into the equations
\begin{equation}
\label{input-output}
\begin{bmatrix}
j_1 \\ 0
\end{bmatrix} = \begin{bmatrix}
D_1 \\ D_2
\end{bmatrix} \nabla G (D_1^\top p_1 + D_2^\top p_2).
\end{equation}
Note that for given $p_1$ and $p_2$, the currents $j_1$ are uniquely defined by the first $N_1$ rows of \eqref{input-output}. In addition, the last $N_2$ rows of \eqref{input-output} provide a relation between $p_1$ and $p_2$. Indeed, 
\begin{equation}
\label{gradientzero}
0 = D_2 \nabla G (D_1^\top p_1 + D_2^\top p_2).
\end{equation}
We will view \eqref{gradientzero} as an input-output model with input $p_1$ and corresponding output $p_2$. The following lemma clarifies the existence and uniqueness of such an output $p_2$, given $p_1$.

\begin{lemma}
\label{l:uniqueoutput}
Assume that the function $G_k : \mathbb{R} \to \mathbb{R}$ is differentiable, coercive and strictly convex for all $k = 1,2,\dots,B$. Then, for every $p_1 \in \R^{N_1}$, the function $H: \mathbb{R}^{N_2} \to \mathbb{R}$ defined as 
\begin{equation*}
H(p_2) = G(D_1^\top p_1 + D_2^\top p_2)
\end{equation*}
is differentiable, coercive and strictly convex. Moreover, for every $p_1 \in \mathbb{R}^{N_1}$, there exists a unique $p_2 \in \mathbb{R}^{N_2}$ satisfying \eqref{gradientzero}. This $p_2$ is the unique minimizer of $H$. 
\end{lemma}

\begin{proof}
Let $p_1 \in \R^{N_1}$. We will first establish the properties of $H$. Obviously, $H$ is differentiable by differentiability of $G$. Moreover, since $\mathcal{G}$ is connected, it follows from \eqref{kernelD} that $D_2^\top$ has full column rank. Indeed, suppose that $D_2^\top x = 0$ for some vector $x \in \R^{N_2}$. Then 
\begin{equation*}
\begin{bmatrix}
0 \\ x
\end{bmatrix} \in \ker D^\top,
\end{equation*}
which means that $x = 0$ by \eqref{kernelD}. Hence, $D_2^\top$ has full column rank. As such, if $\|p_2\| \to \infty$, then also $\|D_1^\top p_1 + D_2^\top p_2\| \to \infty$. This shows that $H$ is coercive since $G$ is coercive. Next, note that for $\alpha \in [0,1]$ and $p_2,\bar{p}_2 \in \mathbb{R}^{N_2}$, we have that $H(\alpha p_2 + (1-\alpha) \bar{p}_2)$ is equal to 
\begin{equation}
\label{sconvH}
\begin{aligned} 
&G(\alpha (D_1^\top p_1 + D_2^\top p_2) + (1-\alpha) (D_1^\top p_1 + D_2^\top \bar{p}_2)) \leq \\ &\alpha G(D_1^\top p_1 + D_2^\top p_2) + (1-\alpha) G(D_1^\top p_1 + D_2^\top \bar{p}_2) = \\
&\alpha H(p_2) + (1-\alpha) H(\bar{p}_2),
\end{aligned}
\end{equation} 
where the inequality follows by convexity of $G$. Thus, also $H$ is a convex function. Moreover, if $\alpha \in (0,1)$ and $p_2 \neq \bar{p}_2$ then the inequality in \eqref{sconvH} is strict due to strict convexity of $G$ and full column rank of $D_2^\top$. This proves that also $H$ is strictly convex. By Proposition~\ref{prop:basicconvexity}~\ref{item:global} in Appendix~\ref{app:convexity}, $H$ has precisely one minimizer. As such, there is precisely one point $p_2$ at which $\nabla H(p_2) = 0$. Now, note that 
\begin{equation*}
\nabla H(p_2) = D_2 \nabla G(D_1^\top p_1 + D_2^\top p_2),
\end{equation*}
which equals the right hand side of \eqref{gradientzero}. We conclude that for every $p_1 \in \mathbb{R}^{N_1}$, there indeed exists a unique $p_2 \in \mathbb{R}^{N_2}$ satisfying \eqref{gradientzero}. 
\end{proof}

\section{The problem of learning in resistive circuits}
\label{sec:problem}

The use of resistive networks in analog computing applications has a long history, see, e.g., \cite{Harris1989}. In this paper, we are interested in \emph{learning} in resistive networks, where we consider the conductances of variable resistors as trainable parameters. We will follow the modeling framework as outlined in the previous section. In particular, we consider a connected undirected graph $\calG = (\calV,\calE)$ with $N$ nodes and $B$ branches. This time, the node set $\calV$ is partitioned as $\calV = \calV_I \cup \calV_H \cup \calV_O$ with $\calV_I$, $\calV_H$ and $\calV_O$ mutually disjoint. We refer to the nodes in $\calV_I$ as \emph{input nodes}, those in $\calV_H$ as \emph{hidden nodes}, and those in $\calV_O$ as \emph{output nodes}. 

All branches of the network contain (nonlinear) \emph{adjustable} resistors. This means that we associate to each of the resistors a \emph{parameterized} co-content function of the form $\theta_k G_k( \cdot )$, where $\theta_k \in \R$ is an adjustable parameter and $G_k : \R \to \R$ for all $k \in \{1,2,\dots,B\}$. We make the following blanket assumption on the functions $G_k$. 

\begin{assumption} 
\label{as:Gk}
The functions $G_k: \R \to \R$ are coercive, twice continuously differentiable, and such that $\frac{d^2 G_k}{dv_k^2}(v_k) > 0$ for all $v_k \in \R$ and $k = 1,2,\dots,B$.
\end{assumption} 

Next, we collect the parameters in the vector
\begin{equation*}
\theta = \begin{bmatrix}
\theta_1 \\ \theta_2 \\ \vdots \\ \theta_{B}
\end{bmatrix} \in \Theta.
\end{equation*}

Here, $\Theta$ denotes the constraint set, assumed to be a nonempty, compact and convex subset of $(0,\infty)^{B}$. As such, we assume that all parameters are positive. This naturally generalizes the case of linear resistors, where the parameters correspond to (positive) tunable conductances. The positivity of the parameters also plays an important role in viewing tunable resistive circuits as input-output models. Indeed, note that Assumption~\ref{as:Gk} implies that the co-content functions $\theta_k G_k( \cdot )$ are differentiable, coercive, and strictly convex\footnote{The assumption on the second derivative of $G_k$ is slightly stronger than strict convexity. The role of this stronger assumption will become clear in Section~\ref{sec:contdif} that discusses the smoothness of certain voltage potentials.} for all positive $\theta_k$. Therefore, Lemma~\ref{l:uniqueoutput} can be readily applied (with $G_k$ replaced by $\theta_k G_k$). 

\begin{remark}
\label{rem:partparameterized}
For notational convenience, we consider all resistors to be parameterized. However, within this framework, we can easily consider the situation in which a subset of resistors is adjustable and the remaining resistors are fixed. This can be done by defining $B_1,B_2 \in \N$ such that $B = B_1 + B_2$, and taking
\begin{equation*}
\Theta := \Theta_1 \times \{1\}^{B_2},
\end{equation*}
where $\Theta_1 \subseteq (0,\infty)^{B_1}$. A relevant setting is the combination of $B_1$ adjustable linear resistors, with parameters correspond to tunable conductances, and $B_2$ fixed nonlinear resistors, see also \cite{kendall2020training}.
\end{remark}

With this in mind, we define the ``parameterized" total co-content function $G: \Theta \times \R^B \to \R$ as
\begin{equation}
\label{totalcocontent}
G(\theta,v) = \sum_{k = 1}^{B} \theta_k G_k(v_k).
\end{equation} 

Let $N_I,N_H$ and $N_O$ denote the number of input, hidden and output nodes, respectively. Then we partition the vector of voltage potentials $p \in \R^N$, the vector of nodal currents $j \in \R^N$ and the incidence matrix $D \in \R^{N \times B}$ as 
\begin{equation*}
p = \begin{bmatrix}
p_I \\ p_H \\ p_O
\end{bmatrix}, \:\: j = \begin{bmatrix}
j_I \\ j_H \\ j_O
\end{bmatrix}, \:\: D = \begin{bmatrix}
D_I \\ D_H \\ D_O
\end{bmatrix},
\end{equation*} 
where $p_I,j_I \in \R^{N_I}$, $p_H,j_H \in \R^{N_H}$, $p_O,j_O \in \R^{N_O}$ and $D_I \in \R^{N_I \times B}$, $D_H \in \R^{N_H \times B}$ and $D_O \in \R^{N_O \times B}$. By \eqref{staticnetwork}, the network model is now given by 
\begin{equation}
\label{networkIHO}
\begin{bmatrix}
j_I \\ j_H \\ j_O
\end{bmatrix} = \begin{bmatrix}
D_I \\ D_H \\ D_O
\end{bmatrix}
\frac{\partial G}{\partial v}(\theta,D^\top_I p_I + D^\top_H p_H + D^\top_O p_O). 
\end{equation}
Setting the nodal current at the hidden and output nodes to zero, i.e., $j_H = 0$ and $j_O = 0$, we have
\begin{equation}
\label{ihomodel}
0 = \begin{bmatrix}
D_H \\ D_O
\end{bmatrix}
\frac{\partial G}{\partial v}(\theta,D^\top_I p_I + D^\top_H p_H + D^\top_O p_O). 
\end{equation}
By Lemma~\ref{l:uniqueoutput}, for every parameter vector $\theta \in \Theta$ and every vector of input potentials $p_I \in \R^{N_I}$, there exists a unique pair of hidden and output potentials $(p_H(\theta,p_I),p_O(\theta,p_I)) \in \R^{N_H} \times \R^{N_O}$ satisfying \eqref{ihomodel}. Here, we stress the dependence of both $p_H$ and $p_O$ on $\theta$ and $p_I$. 

Now, suppose that we have \emph{input-output data} in the form of a single, given pair of potentials 
\begin{equation}
\label{data}
(p_I^d,p_O^d) \in \R^{N_I} \times \R^{N_O}.
\end{equation}
The goal of learning is to fit the parameters to the training data, i.e., we want to compute a parameter vector $\theta^* \in \Theta$ such that 
\begin{equation}
\label{goal}
p_O(\theta^*,p_I^d) \approx p_O^d.
\end{equation}
In other words, the objective is to find a parameter vector $\theta^*$ that \emph{explains} the data in the sense that the network approximately maps $p_I^d$ to the target output $p_O^d$. For simplicity, throughout the paper we focus on a \emph{single} pair of potentials \eqref{data}. This is a simplification of a traditional machine learning setting, that generally assumes \emph{multiple} samples\footnote{In the case of linear resistive networks without hidden nodes, the convergence of contrastive learning was already established for multiple data samples \cite[Sec. V]{huijzer2026convergence}. There, a stochastic version of projected gradient descent was employed. It is expected that a similar approach can be used to generalize the results of this paper to multiple samples; however, this is not the focus of this work.}
$$
(p_I^{d,j},p_O^{d,j}) \in \R^{N_I} \times \R^{N_O}
$$
for $j = 1,2,\dots,n$.

In this work, we study the problem of Eq.~\eqref{goal} through the lens of Target-Clamped Contrastive Learning (TCCL), but before delving into the learning framework (Section~\ref{sec:CL}), we briefly review some fundamental concepts of machine learning and motivate the use of TCCL over Equilibrium Propagation in our setting.

In machine learning, training is usually formalized by introducing a cost function
$$
C : \R^{N_O} \times \R^{N_O} \to \R
$$
with the property that, for all $y \in \R^{N_O}$, $C(\cdot,y)$ has a unique global minimizer, equal to $y$. The cost function $C$ is often chosen to be convex. Examples include the squared error often used for regression tasks, or the categorical cross-entropy used for classification tasks. Given $C$, the objective function $J : \Theta \to \R$ to be minimized with respect to $\theta$ is the average of cost values over the dataset, i.e.,
\begin{equation}
\label{functionJ}
J(\theta) = \frac{1}{n} \sum_{j =1}^n C(p_O(\theta,p_I^{d,j}), p_O^{d,j}).
\end{equation}
To achieve this, the gradient $\nabla J(\theta)$ is usually computed using backpropagation, when applicable. However, in a resistive network, it is unclear how to physically realize backpropagation. For energy-based models like the resistive circuits studied here, alternative approaches have been developed in \cite{kendall2020training,dillavou2022demonstration,anisetti2024frequency,huijzer2026convergence}, sometimes gathered under the umbrella name of \emph{energy-based learning} \citep{Scellier2023}. The main thrust of these works is to replace digital computation of $\nabla J(\theta)$ by learning rules based on \emph{local physical quantities}, thus amenable to efficient implementation in electrical circuits. To achieve this, these algorithms typically compare two network states corresponding to two different boundary conditions. These algorithms either approximate $\nabla J(\theta)$ \citep{kendall2020training,anisetti2024frequency}, compute the exact gradient of an alternative (contrastive) objective function \citep{huijzer2026convergence}, or are based on other heuristics \citep{dillavou2022demonstration}.

In the context of resistive circuits, learning was first formalized using \emph{equilibrium propagation (EP)} \cite{kendall2020training}, a training framework for energy-based models introduced in \cite{scellier2017equilibrium}. EP approximates $\nabla J(\theta)$ by comparing the inference state (called \emph{free} state) and a \emph{nudged} state in which outputs are nudged towards the data. EP requires choosing the nudging strength, and serves as a surrogate for the minimization of $J$.

Contrastive Learning with clamped outputs, which in this work we call \emph{Target-Clamped Contrastive learning} (TCCL), is an alternative energy-based learning rule, which was originally formulated in Hopfield networks \cite{Movellan1991} and has recently been adapted to resistive circuits \cite{huijzer2026convergence}. In TCCL, a \emph{clamped state}, in which outputs are fully clamped to the data (targets), is considered in lieu of the nudged state. This no longer requires a cost function $C$ and a nudging strength; instead, TCCL performs gradient descent on a contrastive function $Q$. Importantly, it was shown in \cite{huijzer2026convergence} that the contrastive function is \emph{convex} for linear resistive circuits without hidden nodes, which is not true for $J$ (see Proposition~\ref{prop:Jnotconvex} of the current paper). The convexity of $Q$ implies that the TCCL algorithm of \cite{huijzer2026convergence} always converges to a \emph{global minimum} of the contrastive function. 

The convexity of the contrastive function in the linear input-output case motivates the study of TCCL for more general resistive networks. The purpose of this paper is to generalize TCCL to nonlinear networks that also include hidden nodes. We will establish a rigorous mathematical framework that will reveal the cases in which i) the contrastive function is convex, and ii) TCCL converges to a global minimum of the contrastive function. In addition, we will study approximation properties of nonlinear resistive networks, showcasing the class of functions that can be approximated by these circuits.

\section{The target-clamped contrastive learning framework}
\label{sec:CL}

In this section, we generalize the contrastive learning framework of \cite{huijzer2026convergence} to the setting of nonlinear resistive networks with input, hidden and output nodes. In Section~\ref{sec:energy}, we will introduce a suitable energy function, and discuss the two energy levels of contrastive learning. In Section~\ref{sec:contdif} we discuss smoothness properties of the voltage potentials arising in these two energy levels. Finally, in Section~\ref{sec:contrastivefunction} we will introduce the contrastive function that will serve as an objective function.

\vspace{5pt}

\subsection{The energy function}
\label{sec:energy}

\noindent
We define the ``energy" function 
\begin{equation*}
E: \Theta \times \mathbb{R}^{N_I} \times \R^{N_H} \times \mathbb{R}^{N_O} \to \mathbb{R}
\end{equation*}
by
\begin{equation}
\label{energyfunction}
E(\theta,p_I,p_H,p_O) := G(\theta,D_I^\top p_I + D_H^\top p_H + D_O^\top p_O).
\end{equation}
Note that the terminology ``energy" is common in the literature on energy-based learning, for example, in the context of Hopfield networks \cite{Movellan1991} and Boltzmann machines \cite{ackley1985learning}. Strictly speaking, however, such functions generally do not correspond to physical energy. In our setting, the function $E$ is equal\footnote{Note that $E$ is a function of $\theta, p_I, p_H$ and $p_O$ while $G$ is a function of $\theta$ and $v$.} to the total co-content of all resistors in the circuit. In the case of linear resistors, $E$ boils down to half the \emph{dissipated power} in the resistors.

Now, consider the data $(p_I^d,p_O^d)$ and a parameter vector $\theta \in \Theta$. As in other models where target-clamped contrastive learning has been employed \citep{Movellan1991,Scellier2023,huijzer2026convergence}, we will consider two energy levels:
\begin{itemize}
\item The \emph{free energy level} $E(\theta,p_I^d,p_H^\ast(\theta),p_O^\ast(\theta))$. Here, $(p_H^\ast(\theta),p_O^\ast(\theta))$ is the pair of so-called \emph{free} hidden and output potentials, defined as the unique minimizer of $E$ for fixed $p_I^d$ and given $\theta$, i.e.,
\begin{equation}
\label{freepotentials}
\:\:\:\:\: (p_H^\ast(\theta),p_O^\ast(\theta)) := \!\!\!\!\!\!\! \argmin_{(x,y) \in \R^{N_H} \times \R^{N_O}} \!\!\! E(\theta,p_I^d,x,y).
\end{equation}
When the dependence on $\theta$ is clear from the context, we will simply write $p_H^\ast$ and $p_O^\ast$. Note that these free potentials satisfy 
\begin{equation*}
0 = \begin{bmatrix} \frac{\partial E}{\partial p_H}(\theta,p_I^d,p_H^\ast,p_O^\ast) \\ \frac{\partial E}{\partial p_O}(\theta,p_I^d,p_H^\ast,p_O^\ast) \end{bmatrix},
\end{equation*}
and thus,
\begin{equation*} 
0 =  
\begin{bmatrix}
D_H \\ D_O
\end{bmatrix}
\frac{\partial G}{\partial v}(\theta,D^\top_I p_I^d + D^\top_H p_H^\ast + D^\top_O p_O^\ast).
\end{equation*}
As such, $p_H^\ast$ and $p_O^\ast$ are precisely the hidden and output potentials arising from the network model \eqref{networkIHO} by setting $j_H = 0$, $j_O = 0$ and $p_I = p_I^d$. In other words, in the notation of Section~\ref{sec:problem}, we have $p_H^\ast(\theta) = p_H(\theta,p_I^d)$ and $p_O^\ast(\theta) = p_O(\theta,p_I^d)$. We recall that the existence and uniqueness of $(p_H^\ast,p_O^\ast)$ follows from Lemma~\ref{l:uniqueoutput} by setting
\begin{equation*}
p_1 = p_I^d \:\: \text{ and } \:\: p_2 = \begin{bmatrix}
p_H^\ast \\ p_O^\ast
\end{bmatrix}.
\end{equation*}
\item The \emph{clamped energy level} $E(\theta,p_I^d,p_H^c(\theta),p_O^d)$, where both input and output potentials are fixed at the data points. Here, $p_H^c(\theta) \in \R^{N_H}$ is the unique minimizer of $E$ for fixed $p_I^d$ and $p_O^d$ and given $\theta$, i.e.,
\begin{equation}
\label{hiddenpotentials}
p_H^c(\theta) := \argmin_{x \in \R^{N_H}} E(\theta,p_I^d,x,p_O^d).
\end{equation}
We will simply write $p_H^{c}$ whenever the dependence on $\theta$ is clear. Note that $p_H^c$ satisfies
\begin{equation*}
0 = \frac{\partial E}{\partial p_H}(\theta,p_I^d,p_H^c,p_O^d),
\end{equation*}
and thus
\begin{equation*}
0 = 
D_H 
\frac{\partial G}{\partial v}(\theta,D^\top_I p_I^d + D^\top_H p_H^c + D^\top_O p_O^d).
\end{equation*}
As such, $p_H^c$ is precisely the vector of hidden potentials arising from the network model \eqref{networkIHO} by setting $j_H = 0$, $p_I = p_I^d$ and $p_O = p_O^d$. Note that the existence and uniqueness of $p_H^c$ again follows from Lemma~\ref{l:uniqueoutput}.
\end{itemize}

\subsection{On the smoothness of the clamped and free potentials}
\label{sec:contdif}

Due to the assumptions on the functions $G_k$, we can say more about the smoothness of the clamped and free potentials, as functions of $\theta$. This will become useful to derive properties of the contrastive function in Section~\ref{sec:contrastivefunction}. Specifically, by Assumption~\ref{as:Gk} the function $(\theta,p_H) \mapsto \frac{\partial E}{\partial p_H}(\theta,p_I^d,p_H,p_O^d)$ is continuously differentiable. Furthermore, the matrix of second derivatives,
\begin{equation*}
\frac{\partial^2 E}{\partial p_H \partial p_H^\top}(\theta,p_I^d,p_H,p_O^d),
\end{equation*} 
is equal to 
\begin{equation*}
D_H 
\frac{\partial^2 G}{\partial v \partial v^\top}(\theta,D_I^\top p_I^d + D_H^\top p_H + D_O^\top p_O^d) D_H^\top.
\end{equation*}
The latter matrix is positive definite for all $\theta \in \Theta$ and $p_H \in \R^{N_H}$. Indeed, this follows from the facts that $D_H$ has full row rank by connectedness of $\mathcal{G}$, and $\frac{\partial^2 G}{\partial v \partial v^\top}(\theta,v)$ is positive definite for all $\theta \in \Theta$ and $v \in \R^B$ by Assumption~\ref{as:Gk}. Now, by the implicit function theorem \cite[Thm. 9.28]{Rudin1976}, $p_H^c(\theta)$ is continuously differentiable and its Jacobian $\frac{\partial p_H^c}{\partial \theta^\top}(\theta)$ is given by 
\begin{equation}
\label{dpHc}
 - \left(D_H \frac{\partial^2 G}{\partial v \partial v^\top}(\theta, v^c(\theta)) D_H^\top \right)^{-1} D_H \frac{\partial^2 G}{\partial v \partial \theta^\top}(\theta, v^c(\theta)),
\end{equation} 
where 
\begin{equation}
\label{defvc}
v^c(\theta) := D^\top \begin{bmatrix}
p_I^d \\ p_H^c(\theta) \\ p_O^d
\end{bmatrix}.
\end{equation}
Similarly, it can be shown that $p_H^\ast(\theta)$ and $p_O^\ast(\theta)$ are also continuously differentiable, and the Jacobian
\begin{equation*}
\begin{bmatrix}
\frac{\partial p_H^\ast}{\partial \theta^\top}(\theta) \\
\frac{\partial p_O^\ast}{\partial \theta^\top}(\theta) 
\end{bmatrix}
\end{equation*}
is given by
\begin{equation}
\label{dpHO}
-\!\!\left(\! D_{HO} \frac{\partial^2 G}{\partial v \partial v^\top}(\theta, v^\ast(\theta)) D_{HO}^\top \!\right)^{\!\!-1}\!\!\!\! D_{HO} \frac{\partial^2 G}{\partial v \partial \theta^\top}(\theta, v^\ast(\theta))
\end{equation} 
where
\begin{equation}
\label{defvast}
D_{HO} := \begin{bmatrix}
D_H \\ D_O
\end{bmatrix}, \text{ and } 
v^\ast(\theta) := D^\top \begin{bmatrix}
p_I^d \\ p_H^\ast(\theta) \\ p_O^\ast(\theta)
\end{bmatrix}.
\end{equation}

\subsection{Contrastive function and its properties}
\label{sec:contrastivefunction}

Now, we will define the \emph{contrastive function}, that compares the free and clamped energy levels. In particular, given the data $(p_I^d,p_O^d)$, we define the function $Q : \Theta \to \mathbb{R}$ as 
\begin{equation}
\label{contrastivefunction}
Q(\theta) := E(\theta,p_I^d,p_H^c(\theta),p_O^d) - E(\theta,p_I^d,p_H^\ast(\theta),p_O^\ast(\theta)).
\end{equation}

The following theorem highlights some of the important features of the function $Q$. 

\begin{theorem}
\label{prop:Q}
The function $Q$ is continuously differentiable. Moreover, $Q$ and its gradient $\nabla Q$ satisfy:
\begin{enumerate}[label=\roman*)]
\item \label{i:nonnegative} $Q(\theta) \geq 0$ for all $\theta \in \Theta$.
\item \label{i:zerooutput} $Q(\theta^d) = 0$ if and only if $\theta^d \in \Theta$ satisfies 
\begin{equation}
\label{thetad}
p_O^\ast(\theta^d) = p_O^d.
\end{equation}
\vspace{-20pt}
\item \label{i:minzero} $\min_{\theta \in \Theta} Q(\theta) = 0$ if and only if there exists $\theta^d \in \Theta$ satisfying \eqref{thetad}.
\item \label{i:gradient} for all $k \in \{1,2,\dots,B\}$,
\begin{equation*}
\left( \nabla Q(\theta) \right)_k = G_k(v^c_k(\theta)) - G_k(v^\ast_k(\theta)),
\end{equation*}
where $v^c(\theta)$ and $v^\ast(\theta)$ are defined in \eqref{defvc} and \eqref{defvast}, respectively.
\item \label{i:QnablaQ} $Q(\theta) = \theta^\top \nabla Q(\theta)$ for all $\theta \in \Theta$. 
\item \label{i:QzeronablaQzero} $Q(\theta) = 0$ if and only if $\nabla Q(\theta) = 0$. 
\item \label{i:Qposhom1} $Q(\lambda \theta) = \lambda Q(\theta)$ for all $\theta \in \Theta$ and $\lambda > 0$ such that $\lambda \theta \in \Theta$.
\end{enumerate}
\end{theorem}

Before turning our attention to proofs, we first discuss some of the consequences of Theorem~\ref{prop:Q}. By~\ref{i:nonnegative}, $Q$ is nonnegative, while \ref{i:minzero} demonstrates that the minimum of $Q$ is zero if and only if there exists a $\theta^d$ that interpolates the data exactly. In the latter case, by \ref{i:zerooutput}, the global minimizers of $Q$ are precisely those parameter vectors that interpolate the data. Moreover, \ref{i:gradient} shows that the gradient $\nabla Q$ can be computed locally, because its $k$-th entry depends only on the local voltages $v_k^c(\theta)$ and $v_k^\ast(\theta)$ over branch $k$. This property, also shared by EP, ensures that contrastive learning can be implemented locally, as explored further in later sections of this paper. Properties~\ref{i:QnablaQ} and \ref{i:QzeronablaQzero} show interesting relations between $Q$ and $\nabla Q$. In particular, under the assumption that the data can be interpolated exactly, \ref{i:QzeronablaQzero} and \ref{i:minzero} show that $\theta^\ast \in \Theta$ is a global minimizer of the \emph{constrained} optimization problem $\min_{\theta \in \Theta} Q(\theta)$ if and only if $\nabla Q(\theta^\ast) = 0$, i.e., the first-order necessary condition for optimality holds for the \emph{unconstrained} optimization problem. Finally, \ref{i:Qposhom1} makes a connection to the concept of \emph{positive homogeneity}. To elaborate on this, we recall the following definition of positive homogeneity, here adapted from \cite{Fleming1977} to functions defined on a domain $\calC \subseteq \R^n$. 

\begin{definition}
Let $k$ be a nonnegative integer and let $\calC \subseteq \R^n$. We call a function $f: \calC \to \R$ \emph{positively homogeneous of degree $k$} if $f(\lambda x) = \lambda^k f(x)$ for all $x \in \calC$ and $\lambda > 0$ such that $\lambda x \in \calC$. 
\end{definition}

The interpretation of Theorem~\ref{prop:Q}~\ref{i:Qposhom1} is thus that $Q: \Theta \to \R$ is \emph{positively homogeneous\footnote{For continuously differentiable functions $f:\R^n \to \R$, the classical \emph{Euler's theorem} (see, e.g., \cite[Thm. 2.1]{Dobbs2018}) asserts that $f$ is positively homogeneous of degree $k$ if and only if $kf(x) = x^\top \nabla f(x)$ for all $x \in \R^n$. Note that Theorem~\ref{prop:Q}~\ref{i:QnablaQ} highlights this characterization (with $k =1$) for the function $Q$ defined on the set $\Theta$.} of degree~$1$}. 

Next, we proceed with the proof of Theorem~\ref{prop:Q}. To this end, we first establish a useful property of the voltage potentials $p_H^\ast, p_O^\ast$ and $p_H^c$. In particular, the next lemma shows that these voltage potentials are positively homogeneous of \emph{degree $0$}. In other words, they are invariant under parameter scaling.

\begin{lemma}
    \label{l:scalingpO}
    The voltage potentials $p_H^\ast: \Theta \to \R^{N_H}$, $p_O^\ast: \Theta \to \R^{N_O}$ and $p_H^c : \Theta \to \R^{N_H}$, as defined in Section~\ref{sec:energy}, are positively homogeneous of degree~$0$. 
\end{lemma}

\begin{proof}
    We only show the proof for $p_H^\ast$ and $p_O^\ast$, as the proof for $p_H^c$ is analogous. Let $\theta \in \Theta$ and $\lambda \theta \in \Theta$ for some scalar $\lambda >0$. Since $E$ in \eqref{energyfunction} is linear in its first argument, we have
    \begin{align*}
        (p_H^\ast(\theta),p_O^\ast(\theta)) 
        &= \argmin_{(x,y) \in \R^{N_H} \times \R^{N_O}} \lambda E(\theta,p_I^d,x,y) \\
        &= \argmin_{(x,y) \in \R^{N_H} \times \R^{N_O}}  E(\lambda \theta,p_I^d,x,y) \\
        & = (p_H^\ast(\lambda \theta),p_O^\ast(\lambda \theta)).
    \end{align*}
    This proves the lemma. 
\end{proof}

\begin{proof}[Proof of Theorem~\ref{prop:Q}] Continuous differentiability of $Q$ immediately follows from Assumption~\ref{as:Gk} and the continuous differentiability of $p_H^c, p_H^\ast$ and $p_O^\ast$, as established in Section~\ref{sec:contdif}. Property~\ref{i:nonnegative} follows by noting that the definitions of $p_H^c$, $p_H^\ast$ and $p_O^\ast$ imply that 
\begin{equation*}
E(\theta,p_I^d,p_H^c(\theta),p_O^d) \geq  E(\theta,p_I^d,p_H^\ast(\theta),p_O^\ast(\theta)),
\end{equation*}
and thus, $Q(\theta) \geq 0$ for all $\theta \in \Theta$. 

We now prove Property~\ref{i:zerooutput}. For the ``if" part, assume that $\theta^d \in \Theta$ satisfies \eqref{thetad}. This implies that
$$
E(\theta^d,p_I^d,p_H^\ast(\theta^d),p_O^\ast(\theta^d)) = E(\theta^d,p_I^d,p_H^\ast(\theta^d),p_O^d),
$$
and, as such, $p_H^\ast(\theta^d)$ is a minimizer of $E(\theta^d,p_I^d, \cdot ,p_O^d)$. By Assumption~\ref{as:Gk} and connectedness of $\mathcal{G}$, $E(\theta^d,p_I^d, \cdot ,p_O^d)$ is continuous, coercive and strictly convex (see Lemma~\ref{l:uniqueoutput}). By Proposition~\ref{prop:basicconvexity}~\ref{item:global}, $E(\theta^d,p_I^d, \cdot ,p_O^d)$ thus has precisely one global minimizer, showing that $p_H^\ast(\theta^d) = p_H^c(\theta^d)$. Therefore, $Q(\theta^d) = 0$. 

To prove the ``only if" part of Property~\ref{i:zerooutput}, assume that $Q(\theta^d) = 0$ for some $\theta^d \in \Theta$. Then
\begin{equation}
\label{eqenergies}
E(\theta^d,p_I^d,p_H^c(\theta^d),p_O^d) = E(\theta^d,p_I^d,p_H^\ast(\theta^d),p_O^\ast(\theta^d)).
\end{equation} 
By Lemma~\ref{l:uniqueoutput} and Proposition~\ref{prop:basicconvexity}~\ref{item:global}, $E(\theta^d,p_I^d, \cdot , \cdot )$ has precisely one minimizer. This minimizer is given by $(p_H^\ast(\theta^d),p_O^\ast(\theta^d))$, and therefore, \eqref{eqenergies} implies that $p_O^\ast(\theta^d) = p_O^d$.  

Next, we prove Property~\ref{i:minzero}. First, note that the minimum in \ref{i:minzero} is well-defined by the extreme value theorem, as $Q$ is continuous and $\Theta$ is compact. The ``if" statement directly follows from Properties~\ref{i:nonnegative} and \ref{i:zerooutput}. Conversely, for the ``only if" part, assume that $\min_{\theta \in \Theta} Q(\theta) = 0$. Then, there exists a $\theta^d \in \Theta$ such that $Q(\theta^d) = 0$. It then follows from Property~\ref{i:zerooutput} that $\theta^d$ satisfies \eqref{thetad}. 

To prove Property~\ref{i:gradient}, first note that $p_H^c,p_H^\ast$ and $p_O^\ast$ are differentiable by the discussion in Section~\ref{sec:contdif}. Then, by definition of $Q$, note that
\begin{equation*}
\begin{aligned} 
\nabla Q(\theta) &= \frac{\partial E}{\partial \theta}(\theta,p_I^d,p_H^c,p_O^d) \\ &+ \left(\frac{\partial p_H^c}{\partial \theta^\top}(\theta)\right)^\top \frac{\partial E}{\partial p_H}(\theta,p_I^d,p_H^c,p_O^d) \\
 &- \frac{\partial E}{\partial \theta}(\theta,p_I^d,p_H^\ast,p_O^\ast) \\ &- \left( \frac{\partial p_H^\ast}{\partial \theta^\top}(\theta) \right)^\top \frac{\partial E}{\partial p_H}(\theta,p_I^d,p_H^\ast,p_O^\ast) \\ &- \left( \frac{\partial p_O^\ast}{\partial \theta^\top}(\theta) \right)^\top \frac{\partial E}{\partial p_O}(\theta,p_I^d,p_H^\ast,p_O^\ast).
\end{aligned}
\end{equation*}
Next, since 
\begin{equation*}
\begin{aligned}
\frac{\partial E}{\partial p_H}(\theta,p_I^d,p_H^c,p_O^d) &= 0, \\
\frac{\partial E}{\partial p_H}(\theta,p_I^d,p_H^\ast,p_O^\ast) &= 0, \text{ and } \\
 \frac{\partial E}{\partial p_O}(\theta,p_I^d,p_H^\ast,p_O^\ast) &= 0,
\end{aligned} 
\end{equation*}
we see that
\begin{equation*}
\nabla Q(\theta) = \frac{\partial E}{\partial \theta}(\theta,p_I^d,p_H^c,p_O^d) - \frac{\partial E}{\partial \theta}(\theta,p_I^d,p_H^\ast,p_O^\ast).
\end{equation*}

Now, Property~\ref{i:gradient} follows from the definitions of the energy function \eqref{energyfunction} and the total co-content \eqref{totalcocontent}.

To prove Property~\ref{i:QnablaQ}, let $\theta \in \Theta$. Note that it follows from \ref{i:gradient} that, 
\begin{equation*}
\theta^\top \nabla Q(\theta) = \sum_{k = 1}^B \theta_kG_k(v_k^c(\theta)) - \sum_{k = 1}^B \theta_k G_k(v_k^\ast(\theta)).
\end{equation*}
Due to the definitions of $G$ in \eqref{totalcocontent},  $E$ in \eqref{energyfunction}, and $Q$ in \eqref{contrastivefunction}, this is precisely equal to $Q(\theta)$. 

To prove Property~\ref{i:QzeronablaQzero}, note that the ``if" statement follows directly from \ref{i:QnablaQ}. To show the ``only if", suppose that $\theta \in \Theta$ is such that $Q(\theta) = 0$. By Property~\ref{i:zerooutput}, $p_O^\ast(\theta) = p_O^d$. Since $p_H^c(\theta)$ is the unique minimizer of $E(\theta,p_I^d, \cdot ,p_O^d)$, this also implies that $p_H^c(\theta) = p_H^\ast(\theta)$. As such, $v^c(\theta) = v^\ast(\theta)$. We conclude from Property~\ref{i:gradient} that $\nabla Q(\theta) = 0$. 

Finally, to prove \ref{i:Qposhom1}, let $\theta \in \Theta$ and $\lambda > 0$ be such that $\lambda \theta \in \Theta$. By \eqref{contrastivefunction}, $Q(\lambda \theta)$ is equal to 
$$
    E(\lambda \theta,p_I^d,p_H^c(\lambda \theta),p_O^d) - E(\lambda \theta,p_I^d,p_H^\ast(\lambda \theta),p_O^\ast(\lambda \theta)),
$$
which is equal to 
$$
    E(\lambda \theta,p_I^d,p_H^c(\theta),p_O^d) - E(\lambda \theta,p_I^d,p_H^\ast(\theta),p_O^\ast(\theta))
$$
due to Lemma~\ref{l:scalingpO}. Finally, since $E$ is linear in the first argument, the latter equals 
$$
    \lambda E(\theta,p_I^d,p_H^c(\theta),p_O^d) - \lambda E( \theta,p_I^d,p_H^\ast(\theta),p_O^\ast(\theta)),
$$
which equals $\lambda Q(\theta)$. This proves the theorem.
\end{proof}

Theorem~\ref{prop:Q} motivates, in a rigorous way, the minimization of the contrastive function $Q$. Indeed, $Q$ is nonnegative and zero precisely at those parameter vectors that interpolate the data. The next sections of this paper are devoted to the analysis of a contrastive learning algorithm used to minimize $Q$. As announced before, Theorem~\ref{prop:Q}~\ref{i:gradient} will become important to ensure the \emph{local} nature of these algorithms. 

In the next section, we will first study the convexity of the function $Q$.

\section{On the convexity of the contrastive function}
\label{sec:convex}

This section is structured as follows. First, in Section~\ref{sec:nohidden}, we will show that $Q$ is convex in the case that there are no hidden nodes (i.e., $\calV_H = \varnothing$). Thereafter, in Section~\ref{sec:nonconvexityJ} we will demonstrate that the function $J$ in \eqref{functionJ}, that EP optimizes, is not convex even in the case with only input-output nodes. This shows an important advantage of using the contrastive function $Q$ in lieu of $J$ in the setting with input and output nodes only. Finally, in Section~\ref{sec:generalcase} we will demonstrate by means of an example that $Q$ is generally not convex if the resistive circuit contains hidden nodes.

\subsection{Convexity of $Q$ without hidden nodes}
\label{sec:nohidden}

In the case of no hidden nodes ($\calV_H = \varnothing$), $p_H$ is void and the energy function boils down to 
\begin{equation}
\label{energynohidden}
E(\theta,p_I,p_O) = G(\theta,D_I^\top p_I + D_O^\top p_O).
\end{equation}
The \emph{free energy} is now given by $E(\theta,p_I^d,p_O^\ast(\theta))$, where $p_O^\ast(\theta)$ is the unique minimizer of $E(\theta,p_I^d, \cdot )$. In turn, the \emph{clamped energy} simply boils down to $E(\theta,p_I^d,p_O^d)$. Therefore, in this setting the contrastive function $Q$ is given by
\begin{equation}
\label{Qnohidden}
Q(\theta) = E(\theta,p_I^d,p_O^d) - E(\theta,p_I^d,p_O^\ast(\theta)).
\end{equation} 
This allows us to establish the following convexity result in Theorem~\ref{t:convexQ}. Note that this is a generalization of \cite[Lem. 5]{huijzer2026convergence} to the nonlinear setting. 

\begin{theorem}
\label{t:convexQ}
If $\calV_H = \varnothing$, then $Q$ is convex. 
\end{theorem}

To prove Theorem~\ref{t:convexQ}, we first state the following lemma. This result is part of Danskin's theorem \cite[p. 173]{Bertsekas2015}, stated here without assuming compactness of $\mathcal{Y}$. For completeness, we include a short proof.

\begin{lemma}
\label{l:concave}
Let $\calX \subseteq \R^n$ be convex and let $\calY \subseteq \R^m$. Consider $f:\calX \times \calY \to \R$. If $f( \cdot ,y)$ is concave for all $y \in \calY$ and $f(x, \cdot )$ has a minimizer for all $x \in \calX$, then the function $g : \calX \to \R$, defined by 
\begin{equation*}
g(x) := \min_{y \in \calY} f(x,y),
\end{equation*}
is concave.
\end{lemma}

\begin{proof}
For every $\alpha \in [0,1]$ and $x_1,x_2 \in \calX$,
\begin{equation*}
\begin{aligned}
&g(\alpha x_1 + (1- \alpha) x_2) = \min_{y \in \calY} f(\alpha x_1 + (1-\alpha)x_2,y) \\
&\geq \min_{y \in \calY} \left( \alpha f(x_1,y) + (1-\alpha)f(x_2,y)\right) \\
& \geq \alpha \min_{y \in \calY} f(x_1,y) + (1-\alpha) \min_{y \in \calY} f(x_2,y) \\ 
&= \alpha g(x_1) + (1-\alpha) g(x_2),
\end{aligned}
\end{equation*}
where the first inequality follows from concavity of $f( \cdot ,y)$. This proves the lemma. 
\end{proof}

\begin{proof}[Proof of Theorem~\ref{t:convexQ}]
In view of the definitions of $E$ in \eqref{energynohidden} and $G$ in \eqref{totalcocontent}, $E(\theta,p_I^d,p_O^d)$ is linear, and thus convex, in $\theta$.  The function $E( \cdot ,p_I^d,p_O)$ is linear, thus concave, for all $p_O \in \R^{N_O}$. Moreover, $E(\theta,p_I^d, \cdot )$ has a (unique) minimizer for all $\theta \in \Theta$ by Assumption~\ref{as:Gk}. As such, by Lemma~\ref{l:concave}, the function
\begin{equation*}
E(\theta,p_I^d,p_O^\ast(\theta)) = \min_{p_O \in \R^{N_O}} E(\theta,p_I^d,p_O)
\end{equation*}
is concave in $\theta$. Thus, $-E(\theta,p_I^d,p_O^\ast(\theta))$ is convex. We conclude from \eqref{Qnohidden} that $Q$ is convex. This proves the theorem. 
\end{proof}

\begin{remark}
\label{rem:convexnonmonotone}
We note that Assumption~\ref{as:Gk} is only used in the proof of Theorem~\ref{t:convexQ} to show that $E(\theta,p_I^d, \cdot )$ has a minimizer for all $\theta \in \Theta$. As such, Theorem~\ref{t:convexQ} is valid even for (non-monotone) resistors with co-content functions $G_k$ that do not satisfy the conditions of Assumption~\ref{as:Gk}, provided that $E(\theta,p_I^d, \cdot )$ still has a minimizer for all $\theta \in \Theta$. The proof of Theorem~\ref{t:convexQ} carries over to this more general case without any changes.
\end{remark}

\subsection{Nonconvexity of $J$}
\label{sec:nonconvexityJ}

In this section, we will show that the objective function $J$ in \eqref{functionJ} is, in general, \emph{not convex}, even for resistive circuits with only input and output nodes. This is different from the contrastive function $Q$, which was shown to always be convex in the setting without hidden nodes (cf. Theorem~\ref{t:convexQ}).

We consider a twice continuously differentiable cost function $C: \R^{N_O} \times \R^{N_O} \to \R$, leading to the objective function $J$ given by
\begin{equation}
\label{functionJ2}
J(\theta) = C(p_O^\ast(\theta),p_O^d).
\end{equation}
Here, we remind the reader that we consider a single data pair $(p_I^d,p_O^d)$, i.e., $n = 1$ in \eqref{functionJ}. In the following proposition, we will work under the assumption that $\Theta$ has nonempty interior (implying that all resistors, including the nonlinear ones, are trainable). In this case, we will show that $J$ can only be convex if it is constant.

\begin{proposition}
\label{prop:Jnotconvex}
    Assume that $\Theta$ has nonempty interior. Suppose that $p_O^\ast:\Theta \to \R^{N_O}$ is twice continuously differentiable. Consider any twice continuously differentiable cost $C: \R^{N_O} \times \R^{N_O} \to \R$. Then, the objective function $J: \Theta \to \R$ in \eqref{functionJ2} is convex if and only if it is constant on $\Theta$.  
\end{proposition}

The proof of Proposition~\ref{prop:Jnotconvex} relies on Lemma~\ref{l:scalingpO}. In fact, we will use the property that $p_O^\ast$ (and thus, $J$) is positively homogeneous of degree $0$.

\begin{proof}[Proof of Proposition~\ref{prop:Jnotconvex}]

    Since the ``if" statement is obvious, we focus on proving the ``only if" part. Assume that $J$ is convex. Then, by the assumptions on $C$ and $p_O^\ast$, the function $J$ is twice continuously differentiable, and 
    $$
    \nabla^2 J(\theta) \geq 0
    $$
    for all $\theta \in \Theta$. Now, let $\theta \in \inter \Theta$ be a point in the interior of $\Theta$. Then, there exists a sufficiently small $\epsilon > 0$ such that $\lambda \theta \in \inter \Theta$ for all $\lambda \in (1-\epsilon,1+\epsilon)$. 
    By Lemma~\ref{l:scalingpO}, $p_O^\ast(\theta) = p_O^\ast(\lambda \theta)$. By definition of $J$, this also implies that 
    \begin{equation}
    \label{poshomJ}
    J(\lambda \theta) = J(\theta)
    \end{equation}
    for all $\lambda \in (1-\epsilon,1+\epsilon)$. 
    
    Differentiating \eqref{poshomJ} on both sides with respect to $\lambda$, and substituting $\lambda = 1$, yields
    \begin{equation}
        \label{thetanablaJ}
    \theta^\top \nabla J(\theta) = 0.
    \end{equation}
    Since $\theta \in \inter \Theta$ was arbitrary, the latter holds for all $\theta \in \inter \Theta$. Now, we differentiate \eqref{thetanablaJ} with respect to $\theta$ to obtain
    \begin{equation}
    \label{derivativesJ}
    \nabla J(\theta) + \nabla^2 J(\theta) \theta = 0
     \end{equation}
    for all $\theta \in \inter \Theta$. Multiplying \eqref{derivativesJ} from left by $\theta^\top$, and substituting \eqref{thetanablaJ}, we get
    $$
    \theta^\top \nabla^2 J(\theta) \theta = 0
    $$
    for all $\theta \in \inter \Theta$. Since $\nabla^2 J(\theta)$ is positive semidefinite, this yields $\nabla^2 J(\theta) \theta = 0$, and therefore
    $$
    \nabla J(\theta) = 0 
    $$
    for all $\theta \in \inter \Theta$, by \eqref{derivativesJ}. As such, $J$ is constant on $\inter \Theta$. By continuity of $J$, the function is also constant on $\Theta$, proving the proposition.
\end{proof}

Under the assumption that $\Theta$ has nonempty interior, Proposition~\ref{prop:Jnotconvex} shows that any reasonable (non-constant) objective function $J$ is \emph{not} convex, even if the cost function $C$ is convex. In case $\inter \Theta$ is empty
(which happens for example if some resistors are fixed and untrainable)
, $J$ may be convex or not, depending on the particular constraint set $\Theta$, the graph structure, co-content functions and the data. As mentioned before, Proposition~\ref{prop:Jnotconvex} highlights an advantage of the contrastive function $Q$ as objective in lieu of the objective $J$ that Equilibrium Propagation optimizes. 

\subsection{Nonconvexity of $Q$ in the general case}
\label{sec:generalcase}

It is quite straightforward to demonstrate that $Q$ is, in general, not convex in the presence of hidden nodes. In the following, we provide a concrete example.

\begin{example}
\label{ex:nonconvex}
Consider the circuit with two input, two hidden, and two output nodes, given in Figure~\ref{fig:nonconvex}. We assume all resistors are linear (so $G_k(v_k) = \frac{1}{2}v_k^2$ for $k = 1,2,\dots,6$), and consider the data
\begin{equation*}
p_I^d = \begin{bmatrix}
1 \\ -3
\end{bmatrix}, \quad p_O^d = \begin{bmatrix}
-2.0968 \\
   -2.3548
\end{bmatrix}.
\end{equation*}
In this case, $\nabla^2 Q(\ones)$ has a negative eigenvalue, showing that $Q$ is not convex. 
\begin{figure}[h!]
\centering
\includegraphics[width=0.4\textwidth]{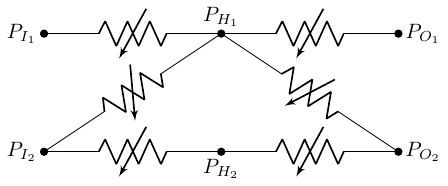}
\caption{Example of circuit with non-convex $Q$ used in Example~\ref{ex:nonconvex}.}
\label{fig:nonconvex}
\end{figure}
\end{example}

\section{The target-clamped contrastive learning algorithm}
\label{sec:algorithm}

In what follows, we will assume that 
\begin{equation}
\label{defTheta}
\Theta = \{ x \in \R^{B} \mid \delta_i \leq x_i \leq \epsilon_i \text{ for all } i = 1,2,\dots,B \},
\end{equation} 
where $\delta_i, \epsilon_i \in \R$ are given constants satisfying the inequalities $0 < \delta_i \leq \epsilon_i$ for all $i = 1,2,\dots,B$. 

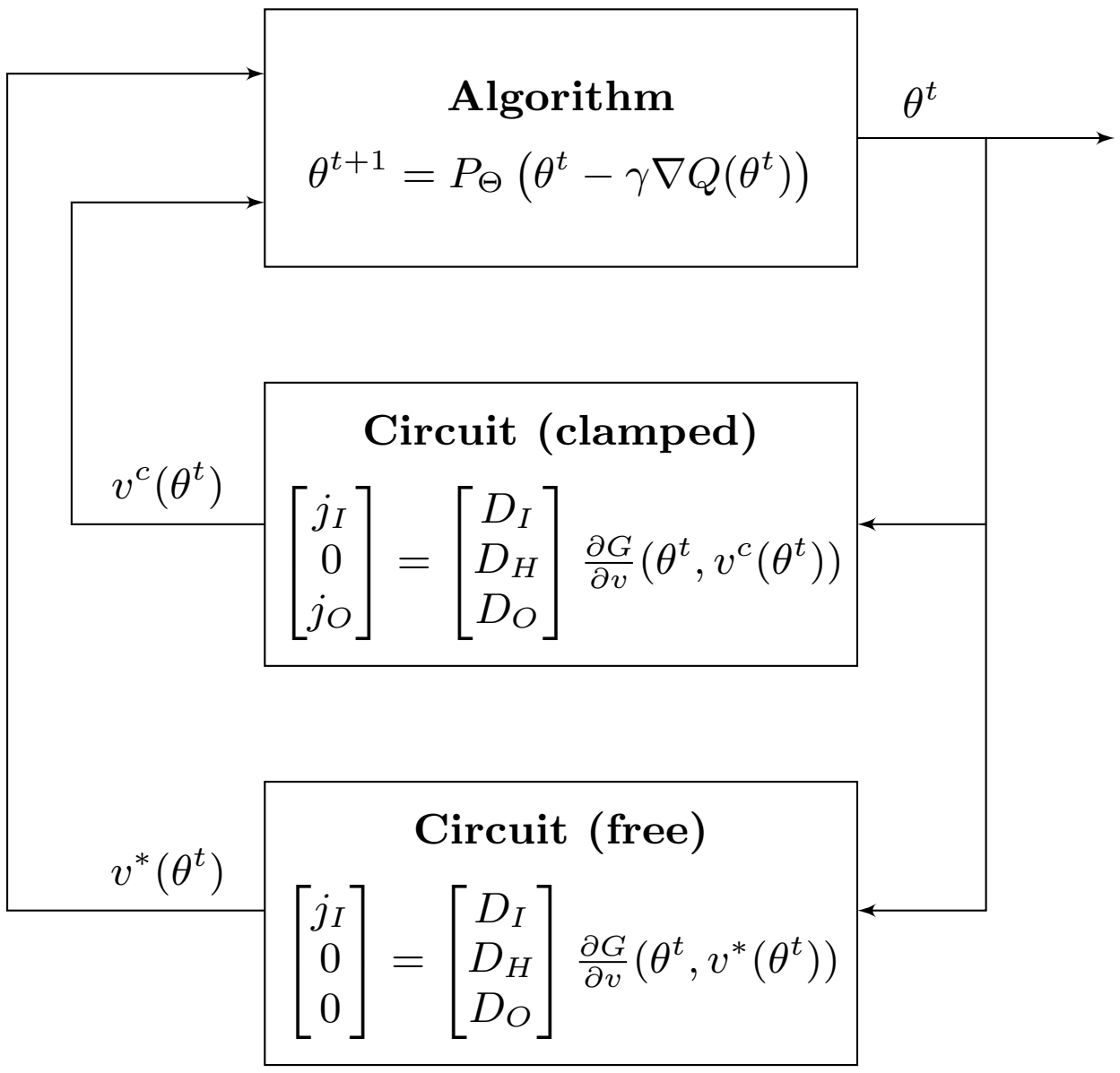
\begin{figure}[h!]
        \begin{center}  
            \begin{tikzpicture}[auto, node distance=0.5cm,>=latex',scale=0.95]
            \draw (-0.8,0) rectangle (3.8,-2);
            \node at (1.5,-0.7) {{\bf Algorithm}};
            \node at (1.5,-1.3) {$\theta^{t+1} = P_{\Theta} \left( \theta^t - \gamma \nabla Q(\theta^t) \right)$};
			\draw (-0.8,-2.9) rectangle (3.8,-5.1);
			\node at (1.5,-3.3) {{\bf Circuit (clamped)}};
			\node at (1.5,-4.3) {$\begin{bmatrix}
j_I \\ 0 \\ j_O
\end{bmatrix} = \begin{bmatrix}
D_I \\ D_H \\ D_O
\end{bmatrix}
\frac{\partial G}{\partial v}(\theta^t,v^c(\theta^t))$};
			\draw (-0.8,-6) rectangle (3.8,-8.2);
			\node at (1.5,-6.4) {{\bf Circuit (free)}};
			\node at (1.5,-7.4) {$\begin{bmatrix}
j_I \\ 0 \\ 0
\end{bmatrix} = \begin{bmatrix}
D_I \\ D_H \\ D_O
\end{bmatrix}
\frac{\partial G}{\partial v}(\theta^t,v^\ast(\theta^t))$};
			\draw[->] (3.8,-1) -- (4.8,-1) -- (4.8,-4) -- (3.8,-4);
			\draw[->] (4.8,-1) -- (5.8,-1);
			\draw[->] (4.8,-4) -- (4.8,-7) -- (3.8,-7);
			\draw[->] (-0.8,-4) -- (-2.3,-4) -- (-2.3,-1.5) -- (-0.8,-1.5);
			\draw[->] (-0.8,-7) -- (-2.8,-7) -- (-2.8,-0.5) -- (-0.8,-0.5);
			\node at (4.3,-0.7) {$\theta^t$};
			\node at (-1.55,-3.7) {$v^c(\theta^t)$};
			\node at (-1.55,-6.7) {$v^\ast(\theta^t)$};
            \end{tikzpicture}
        \end{center}
        \caption{Illustration of the proposed target-clamped contrastive learning (TCCL) algorithm and its interaction with the circuit. To perform an update of the parameters at step $t$ of the algorithm, we need two different ``experiments", resulting in the voltages $v^c(\theta^t)$ and $v^\ast(\theta^t)$. In both experiments, the $k$-th entry of these voltages can be measured locally at branch $k$ of the circuit, enabling the local update rule \eqref{localupdate}.}
        \label{fig:algorithm}
\end{figure}

The projection onto $\Theta$ takes the form 
\begin{equation}
\label{expressionP}
(P_\Theta(x))_k = \begin{cases} x_k & \text{ if } \delta_k \leq x_k \leq \epsilon_k \\ \delta_k & \text{ if } x_k < \delta_k \\ \epsilon_k & \text{ otherwise}, \end{cases}
\end{equation} 
for all $k = 1,2,\dots,B$. We now consider the following \emph{projected gradient descent} algorithm:
\begin{equation}
\label{EBLalgorithm}
\theta^0 \in \Theta \:\: \text{ and } \:\: \theta^{t+1} = P_{\Theta} \left( \theta^t - \gamma \nabla Q(\theta^t) \right) \:\:\text{ for } \:\: t \in \mathbb{N}.
\end{equation}
Here $\gamma > 0$ is the to-be-determined stepsize of the algorithm. This algorithm was already studied in \cite{huijzer2026convergence} for linear resistive circuits without hidden nodes. Here we treat the nonlinear case, including hidden nodes. A key feature of \eqref{EBLalgorithm} is that it is a contrastive learning algorithm that can be implemented \emph{locally} due to Theorem~\ref{prop:Q}~\ref{i:gradient} and \eqref{expressionP}. Specifically, \eqref{EBLalgorithm} boils down to $\theta_k^0 \in [\delta_k,\epsilon_k]$ and 
\begin{equation}
\label{localupdate}
\theta_k^{t+1} = P_{[\delta_k,\epsilon_k]} \left( \theta_k^t - \gamma (G_k(v_k^c(\theta^t)) - G_k(v_k^\ast(\theta^t))) \right)
\end{equation} 
for $t \in \N$ and $k = 1,2,\dots,B$. Note that the update of the parameter $\theta_k^t$ requires two (local) voltages $v_k^c(\theta^t)$ and $v_k^\ast(\theta^t)$ at branch $k$. At each step $t$, these voltages can be obtained by two different ``experiments". The first one clamps (i.e., fixes) the input and output potentials at $p_I^d$ and $p_O^d$ and takes $j_H = 0$, resulting in the clamped voltage $v^c_k(\theta^t)$. The second one only fixes the input potential $p_I^d$, and takes $j_H = 0$ and $j_O = 0$, resulting in the free voltage $v_k^\ast(\theta^t)$ at branch $k$. The algorithm is further illustrated in Figure~\ref{fig:algorithm}.

The expression \eqref{localupdate} highlights an important feature of the contrastive learning rule, previously also discussed in the context of equilibrium propagation, see e.g. \cite[Sec. II]{Scellier2026quantum}. Specifically, the update rule for parameter $\theta_k$ depends explictly on the co-content function $G_k$ at branch $k$. Its dependence on the other co-content functions in the network is only implicit, through the voltages $v_k^c$ and $v_k^*$. This is especially relevant for networks with adjustable linear resistors and nonlinear fixed resistors, as also discussed in Remark~\ref{rem:partparameterized}. In this case, the update rule for each of the linear resistors takes the form 
$$
\theta_k^{t+1} = P_{[\delta_k,\epsilon_k]} \left( \theta_k^t - \frac{\gamma}{2} \left( (v_k^c(\theta^t))^2 - (v_k^\ast(\theta^t))^2 \right) \right),
$$
which does not depend on a model of the nonlinear resistors. This highlights a robustness property of contrastive learning: it is possible to implement exact projected gradient descent on the contrastive function, even with incomplete knowledge of the current-voltage characteristics of the nonlinear resistors.

In what follows, we will formally study the convergence of the contrastive learning algorithm \eqref{EBLalgorithm}. To do so, we will first provide conditions under which $\nabla Q$ is Lipschitz continuous in Section~\ref{sec:Lipschitzcont}. Thereafter, we will establish the convergence of contrastive learning without hidden nodes in Section~\ref{sec:convergenceconvex}, relying on the convexity result in Theorem~\ref{t:convexQ}. In Section~\ref{sec:convergencenonconvex}, we will establish a proof of convergence in the more challenging situation including hidden nodes.

\subsection{Lipschitz continuity of $\nabla Q$}
\label{sec:Lipschitzcont}

In this section, we will prove that $\nabla Q$ is Lipschitz continuous, which will be used to establish convergence guarantees of the contrastive learning algorithm.

\begin{proposition}
\label{p:Lipschitzcont}
The Hessian of $Q$ satisfies
\begin{equation}
\label{HesQ}
\begin{aligned} 
\nabla^2 Q(\theta) &= \begin{bmatrix}
\frac{dG_1}{dv_1}(v_1^c(\theta)) & \cdots & 0 \\
\vdots & \ddots & \vdots \\
0 & \cdots & \frac{dG_B}{dv_B}(v_B^c(\theta))
\end{bmatrix} \frac{\partial v^c}{\partial \theta^\top}(\theta) \\ &- \begin{bmatrix}
\frac{dG_1}{dv_1}(v_1^\ast(\theta)) & \cdots & 0 \\
\vdots & \ddots & \vdots \\
0 & \cdots & \frac{dG_B}{dv_B}(v_B^\ast(\theta))
\end{bmatrix} \frac{\partial v^\ast}{\partial \theta^\top}(\theta).
\end{aligned} 
\end{equation}
Moreover, the function $\nabla Q$ is Lipschitz continuous with constant 
\begin{equation}
\label{defL}
L := \max_{\theta \in \Theta} \norm{ \nabla^2 Q(\theta) }.
\end{equation}
\end{proposition}

\begin{proof}
We will first establish the formula for the Hessian of $Q$. Thereafter, the Lipschitz continuity of $\nabla Q$ will follow by \cite[Lem. 3.1]{Khalil2014}. 

We first note that, due to Theorem~\ref{prop:Q}~\ref{i:gradient}, 
\begin{equation*}
\begin{aligned}
\frac{\partial}{\partial \theta_j} (\nabla Q(\theta))_k &= 
\frac{dG_k}{dv_k}(v_k^c(\theta)) \frac{\partial v_k^c}{\partial \theta_j}(\theta) \\ &- \frac{dG_k}{dv_k}(v_k^\ast(\theta)) \frac{\partial v_k^\ast}{\partial \theta_j}(\theta),
\end{aligned}
\end{equation*}
for all $k,j = 1,2,\dots,B$. This means that the Hessian of $Q$ takes the form \eqref{HesQ}. By Assumption~\ref{as:Gk}, $\frac{dG_i}{d v_i}$ is continuous for all $i = 1,2,\dots,B$. Moreover, due to the discussion in Section~\ref{sec:contdif}, $v^c$, $v^\ast$, $\frac{\partial v^c}{\partial \theta^\top}$ and $\frac{\partial v^\ast}{\partial \theta^\top}$ are continuous. As such, $\nabla^2 Q$ is also continuous. Obviously, $\Theta$ is compact due to \eqref{defTheta}. Hence, by the extreme value theorem, the function $\| \nabla^2 Q(\theta) \|$ attains its supremum on $\Theta$. As such, $L$ in \eqref{defL} is well-defined. Finally, because $\| \nabla^2 Q(\theta) \| \leq L$ for all $\theta \in \Theta$, we conclude by \cite[Lem. 3.1]{Khalil2014} that $\nabla Q$ is Lipschitz continuous with constant $L$. 
\end{proof}

\begin{remark} 
In the case of \emph{linear} resistors, i.e., $G_k(v_k) = \frac{1}{2} v_k^2$ for all $k = 1,2,\dots,B$, the formula for the Hessian of $Q$ simplifies. Indeed, in this case, $\frac{dG_k}{dv_k}(v_k) = v_k$ for all $k = 1,2,\dots,B$, while
\begin{equation*}
\frac{\partial^2 G}{\partial v \partial v^\top}(\theta,v) = \diag(\theta), \:\: \text{ and } \:\: \frac{\partial^2 G}{\partial v \partial \theta^\top}(\theta,v) = \diag(v).
\end{equation*}
Thus, by \eqref{dpHc}, \eqref{defvc}, \eqref{dpHO}, and \eqref{defvast}, the Hessian in \eqref{HesQ} simplifies to
\begin{equation*}
\begin{aligned} 
\nabla^2 Q(\theta) &= \diag(v^\ast(\theta)) W_{HO}(\theta) \diag(v^\ast(\theta)) \\ &- \diag(v^c(\theta)) W_H \diag(v^c(\theta)),
\end{aligned} 
\end{equation*} 
where
\begin{equation*}
\begin{aligned} 
W_{HO}(\theta) \! &:= \! \begin{bmatrix}
D_H \\ D_O
\end{bmatrix}^\top\!\! \left( \begin{bmatrix}
D_H \\ D_O
\end{bmatrix} \diag(\theta) \begin{bmatrix}
D_H \\ D_O
\end{bmatrix}^\top \right)^{\! -1} \!\!\begin{bmatrix}
D_H \\ D_O
\end{bmatrix}, \\ W_{H}(\theta) &:= D_H^\top \left( D_H \diag(\theta) D_H^\top \right)^{-1} D_H.
\end{aligned} 
\end{equation*} 
In the simple case that there are no hidden nodes ($\calV_H = \varnothing$), the matrix $D_H$ is void and $\nabla^2 Q(\theta)$ equals 
\begin{equation*}
\diag(v^\ast(\theta)) D_O^\top \left( D_O \diag(\theta) D_O^\top \right)^{-1} D_O \diag(v^\ast(\theta)).
\end{equation*}
Note that the latter matrix is positive semidefinite for all $\theta \in \Theta$. This recovers the argument in \cite{huijzer2026convergence}, showing that $Q$ is convex in the special case of linear circuits without hidden nodes, see \cite[Lem. 5]{huijzer2026convergence}. 
\end{remark}

\subsection{Convergence of contrastive learning without hidden nodes}
\label{sec:convergenceconvex}

In this section, we will establish a basic convergence result for contrastive learning in the case that there are no hidden nodes (i.e., $\calV_H = \varnothing$). For this, we will rely on the convexity result in Theorem~\ref{t:convexQ}. We first recall the following result from the literature, on the convergence of projected gradient descent.

\begin{proposition}
\label{p:PGDconvex}
Let $\calC \subseteq \R^n$ be a nonempty closed convex set, and let $f: \R^n \to \R$ be differentiable and convex. Suppose that $\argmin_{x\in\calC} f(x) \neq \varnothing$ and $\nabla f$ is Lipschitz continuous with constant $L$. Consider the projected gradient descent algorithm
\begin{equation*}
x^0 \in \calC, \:\: \text{ and } \:\: x^{t+1} = P_\calC\left( x^t - \gamma \nabla f(x^t) \right) \:\: \text{ for } t \in \N.
\end{equation*}
If $\gamma \in (0,\frac{2}{L})$ then 
\begin{equation*}
\lim_{t \to \infty} x^t = x^*,
\end{equation*}
where $x^* \in \argmin_{x\in\calC} f(x)$. 
\end{proposition}

Proposition~\ref{p:PGDconvex} is a classical result in convex optimization. Typically, this result is stated for the range of stepsizes $(0,\frac{1}{L}]$, see e.g. \cite[Thm. 10.24]{Beck2017}. One way of showing that projected gradient descent converges also for the larger interval of stepsizes $(0,\frac{2}{L})$ is to apply methods from monotone operator theory, in particular, the Krasnosel'ski\u{\i}-Mann theorem, see e.g. \cite[Thm. 1]{Ryu2022} and \cite[p. 49]{Ryu2022}. For the sake of completeness, we have included a self-contained proof of Proposition~\ref{p:PGDconvex} in Appendix~\ref{sec:proofpropconv}.

Based on Proposition~\ref{p:PGDconvex} we now prove that the contrastive learning algorithm in \eqref{EBLalgorithm} always converges to a global minimizer of the contrastive function, assuming there are no hidden nodes. This is the topic of the next theorem.

\begin{theorem}
\label{t:convergenceconvex}
Assume that $\calV_H = \varnothing$. Define $L$ as in \eqref{defL}. Consider the contrastive learning algorithm in \eqref{EBLalgorithm}. If $\gamma \in (0,\frac{2}{L})$, then, for every $\theta^0 \in \Theta$, 
\begin{equation*}
\lim_{t \to \infty} \theta^t = \theta^*,
\end{equation*}
for some parameter vector $\theta^* \in \argmin_{\theta \in \Theta} Q(\theta)$. Furthermore, if the data can be interpolated exactly, i.e., if \eqref{thetad} holds for some $\theta^d \in \Theta$, then also $p_O^\ast(\theta^\ast) = p_O^d$.
\end{theorem}

\begin{proof}
The set $\Theta$ in \eqref{defTheta} is nonempty, closed and convex. Moreover, the contrastive function $Q$ is differentiable and, since $\calV_H = \varnothing$, convex by Theorem~\ref{t:convexQ}. The function $Q$ has at least one minimizer by the extreme value theorem, as $Q$ is continuous and $\Theta$ is compact. Moreover, the gradient $\nabla Q$ is Lipschitz continuous with constant $L$ given in \eqref{defL}, by Proposition~\ref{p:Lipschitzcont}. Therefore, Theorem~\ref{t:convergenceconvex} now follows from Proposition~\ref{p:PGDconvex}, and Theorem~\ref{prop:Q}~\ref{i:minzero} and \ref{i:zerooutput}.
\end{proof}

\begin{remark}
The convergence rate of the basic projected gradient descent algorithm \eqref{EBLalgorithm} can be improved by considering momentum terms. In \cite{Nguyen2026}, we have provided a convergence guarantee, as the one in Theorem~\ref{t:convergenceconvex}, for contrastive learning with \emph{Nesterov acceleration}. For a fixed step size $\gamma \in (0,\frac{1}{L}]$, the convergence rates of the non-accelerated and Nesterov variants of the algorithm are $\mathcal{O} (\frac{1}{t})$ and $\mathcal{O} (\frac{1}{t^2})$, respectively.
\end{remark}

\subsection{Convergence of contrastive learning with hidden nodes}
\label{sec:convergencenonconvex}

In this section we study the convergence of target-clamped contrastive learning in the more general situation that $\calV_H \neq \varnothing$. As we have seen, in this situation the contrastive function is, in general, not convex. For non-convex optimization problems, some care has to be taken with the convergence of the iterates $\theta^t$, and their convergence does not follow from the convergence of $Q(\theta^t)$, in general \cite{Bertsekas2016}. The following theorem asserts that contrastive learning, applied to resistive circuits with hidden nodes, converges on a subsequence to a stationary point $\theta^\ast$ of the constrained optimization problem $\min_{\theta \in \Theta} Q(\theta)$. Furthermore, interestingly, any such stationary points in the \emph{interior} of $\Theta$ are global minimizers of $Q$, and are such that $p_O^\ast(\theta^\ast) = p_O^d$. 

\begin{theorem}
\label{t:convergencehidden}
Define $L$ as in \eqref{defL}. Consider the contrastive learning algorithm in \eqref{EBLalgorithm}. Let $\theta^0 \in \Theta$ and $\gamma \in (0,\frac{2}{L})$. Then $\lim_{t \to \infty} Q(\theta^t)$ exists. Furthermore, the sequence $(\theta^t)_{t=0}^\infty$ has a convergent subsequence with limit $\theta^\ast \in \Theta$ satisfying
\begin{equation*}
\nabla Q(\theta^\ast)^\top (\theta - \theta^\ast) \geq 0
\end{equation*}
for all $\theta \in \Theta$. If $\theta^* \in \inter \Theta$ then $Q(\theta^*) = 0$ and we have that $p_O^\ast(\theta^\ast) = p_O^d$. 
\end{theorem}

\begin{proof}
Let $t \in \N$. We start by applying Lemma~\ref{l:projection}~\ref{projineq} in Appendix~\ref{app:convexity} with $\calC = \Theta$, $x = \theta^t - \gamma \nabla Q(\theta^t)$ and $y = \theta^t$, which yields
\begin{equation*}
(\theta^t - \gamma \nabla Q(\theta^t) - \theta^{t+1})^\top (\theta^t-\theta^{t+1}) \leq 0.
\end{equation*}
In other words,
\begin{equation}
\label{thetanabla}
(\theta^{t+1} - \theta^t)^\top \nabla Q(\theta^t) \leq -\frac{1}{\gamma} \| \theta^{t+1} - \theta^t \|^2.
\end{equation}
Now, since $\nabla Q$ is Lipschitz continuous with constant $L$ in \eqref{defL}, we obtain, by Lemma~\ref{l:Lipschitzineq} in Appendix~\ref{app:convexity},
\begin{equation}
\label{descentQ}
\begin{aligned}
Q(\theta^{t+1}) - Q(\theta^t) &\leq (\theta^{t+1} - \theta^t)^\top \nabla Q(\theta^t) \\ + \frac{L}{2} \|\theta^{t+1} - \theta^t \|^2
&\leq \left(\frac{L}{2} - \frac{1}{\gamma}\right) \|\theta^{t+1} - \theta^t \|^2,
\end{aligned}
\end{equation}
where the second inequality follows by \eqref{thetanabla}. Since $\gamma \in (0,\frac{2}{L})$, we have $\frac{L}{2} - \frac{1}{\gamma} < 0$. As such, the sequence $(Q(\theta^t))_{t=0}^\infty$ is nonincreasing. By Theorem~\ref{prop:Q}~\ref{i:nonnegative}, $Q$ is nonnegative and thus bounded from below. Hence, by the monotone convergence theorem, $\lim_{t \to \infty} Q(\theta^t)$ exists. 

Next, since $\Theta$ is compact by \eqref{defTheta}, it follows from the Bolzano-Weierstrass theorem that $(\theta^t)_{t=0}^\infty$ has a convergent subsequence $(\theta^{t_j})_{j=0}^\infty$, where $0 \leq t_0 < t_1 < \cdots$. Denote the limit of this subsequence by $\theta^\ast \in \Theta$. By \eqref{descentQ}, $\|\theta^{t+1} - \theta^t\| \to 0$ as $t \to \infty$. Thus, $\theta^{t_j+1} \to \theta^{\ast}$ as $j \to \infty$. Now, take the limit on both sides of 
\begin{equation*}
\theta^{t_j+1} = P_\Theta(\theta^{t_j} - \gamma \nabla Q(\theta^{t_j}))
\end{equation*}
as $j \to \infty$, and use continuity of both $P_\Theta$ and $\nabla Q$ to obtain 
\begin{equation*}
\theta^\ast = P_\Theta(\theta^\ast - \gamma \nabla Q(\theta^\ast)).
\end{equation*}

Now, by Lemma~\ref{l:projection}~\ref{projineq} in Appendix~\ref{app:convexity}, 
\begin{equation}
\label{statpointproof}
\nabla Q(\theta^\ast)^\top (\theta - \theta^\ast) \geq 0
\end{equation} 
for all $\theta \in \Theta$, showing that $\theta^\ast$ is indeed a stationary point.

Finally, suppose that $\theta^\ast \in \inter \Theta$. Then, for sufficiently small $\epsilon > 0$, 
\begin{equation*}
\theta^\ast - \epsilon \nabla Q(\theta^\ast) \in \Theta.
\end{equation*}
Due to \eqref{statpointproof}, 
\begin{equation*}
-\epsilon \|\nabla Q(\theta^\ast)\|^2 \geq 0,
\end{equation*}
meaning that $\nabla Q(\theta^\ast) = 0$. Now, by Theorem~\ref{prop:Q}~\ref{i:QzeronablaQzero}, $Q(\theta^\ast) = 0$. 
Finally, by Theorem~\ref{prop:Q}~\ref{i:zerooutput}, $p_O^\ast(\theta^\ast) = p_O^d$. This proves the theorem. 
\end{proof}

By Theorem~\ref{t:convergencehidden}, any limit point of \eqref{EBLalgorithm} in the interior of $\Theta$ is a \emph{global} minimizer of $Q$. Furthermore, any such point interpolates the data exactly. This also shows that the contrastive learning algorithm can only converge to interior points of $\Theta$ if the data can be interpolated exactly, i.e. if there exists $\theta^d$ such that \eqref{thetad} holds. If not, the algorithm always converges, subsequentially, to stationary points on the boundary of $\Theta$.

\section{Computational expressivity of nonlinear resistive networks with monotone elements}
\label{sec:approximation}

Universal function approximation is a classical topic in the study of artificial neural networks, see e.g. \cite{Cybenko1989}. So far, however, relatively few results have been obtained on the universal approximation capabilities of \emph{analog} circuits \cite{Scellier2025}. In the context of this paper, it is therefore relevant to establish the class of functions that the considered nonlinear resistive circuits with monotone elements can approximate. To this end, we first introduce the setting. 

We consider a class of resistive circuits with a fixed number of input and output nodes, but with varying number of branches and hidden nodes. Specifically, let $N_I, N_O \in \N$ be fixed. We consider the class of resistive circuits defined by 
\begin{enumerate}[label=(\alph*)]
\item \label{i:graph} a connected, undirected graph $\mathcal{G} = (\mathcal{V},\mathcal{E})$, where $\calV = \calV_I \cup \calV_H \cup \calV_O$ with $\calV_I, \calV_H, \calV_O$ mutually disjoint, and $|\calV_I| = N_I$ and $|\calV_O| = N_O$.
\item \label{i:func} a set of co-content functions $G_k : \R \to \R$, where, for all $k = 1,2,\dots,|\calE|$, $G_k$ is 
\begin{itemize}
\item twice continuously differentiable,
\item coercive, and such that
\item $\frac{d^2 G_k}{d v_k^2}(v_k) > 0$ for all $v_k \in \R$.
\end{itemize}
\end{enumerate}
Consider any such graph $\calG = (\calV,\calE)$ and co-content functions $G_k$ ($k=1,2,\dots,|\calE|$). Define 
\begin{equation*}
G(v) = \sum_{k = 1}^{|\calE|} G_k(v_k),
\end{equation*}
for $v \in \R^{|\calE|}$, and let 
\begin{equation*}
D = \begin{bmatrix}
D_I \\ D_H \\ D_O
\end{bmatrix}
\end{equation*}
be the incidence matrix of $\calG$, where $D \in \R^{|\calV| \times |\calE|}$, $D_I \in \R^{N_I \times |\calE|}$, $D_H \in \R^{|\calV_H| \times |\calE|}$ and $D_O \in \R^{N_O \times |\calE|}$. We denote the input-output function of the circuit, mapping a vector of input potentials to its corresponding vector of output potentials, by $F(\cdot)$. In other words, $F$ is the function such that, for all $p_I \in \R^{N_I}$, 
\begin{equation}
\label{networkmodeluni}
\begin{bmatrix}
D_H \\ D_O
\end{bmatrix}\nabla G(D_I^\top p_I + D_H^\top p_H + D_O^\top F(p_I)) = 0
\end{equation}
for some vector of hidden potentials $p_H \in \R^{|\calV_H|}$. Note that $F$ is well-defined by Lemma~\ref{l:uniqueoutput}. 

Now, let 
\begin{equation*}
\calF_{N_I,N_O} \subseteq (\R^{N_O})^{\R^{N_I}}
\end{equation*}
denote the set of all functions from $\R^{N_I}$ to $\R^{N_O}$ that can be obtained this way, i.e., all input-output maps $F$ corresponding to a circuit whose graph satisfies \ref{i:graph} and co-content functions satisfy \ref{i:func}. 

The following lemma establishes two important properties of the functions $F \in \calF_{N_I,N_O}$. Here, we will use the notation $x \preceq y$ for two vectors $x,y \in \R^n$, indicating that $x_i \leq y_i$ for all $i = 1,2,\dots,n$. 

\begin{lemma}
\label{l:propertiesF}
Let $F \in \calF_{N_I,N_O}$. Then
\begin{enumerate}[label=\roman*)]
\item \label{i:xones} for all $x \in \R^{N_I}$ and $\alpha \in \R$, 
\begin{equation*}
F(x + \alpha \ones) = F(x) + \alpha \ones.
\end{equation*}
\item \label{i:xy} for all $x,y \in \R^{N_I}$, 
\begin{equation*}
x \preceq y \implies F(x) \preceq F(y).
\end{equation*}
\end{enumerate}
\end{lemma}

The proof of Lemma~\ref{l:propertiesF} is given in Appendix~\ref{sec:prooflempropF}. Lemma~\ref{l:propertiesF} clearly demonstrates that the considered nonlinear resistive circuits are not universal approximators, in the sense that they cannot approximate any continuous function on a compact set up to an arbitrary precision. The physics of the circuit and the assumption of strictly monotone current-voltage characteristics (Assumption~\ref{as:Gk}) impose the two important constraints \ref{i:xones} and \ref{i:xy} that any input-output map $F$ must exhibit. Property~\ref{i:xones} is related to the fact that adding the same constant to all potentials does not change the network voltages, while Property~\ref{i:xy} expresses a form of monotonicity: if the inputs $x$ and $y$ are ordered as $x \preceq y$, then also the outputs satisfy $F(x) \preceq F(y)$. 

\begin{remark}
The nonlinear resistive networks of \cite{kendall2020training} circumvent the monotonicity constraint imposed by Lemma~\ref{l:propertiesF} by employing a dual-rail representation -- doubling input nodes and inverting one set of voltages. Alternatively, \cite{Scellier2025} implement this dual-rail representation using voltage-controlled voltage sources (VCVSs), and establish that these nonlinear resistive networks with monotone elements and VCVSs are universal function approximators.
\end{remark}

In view of Lemma~\ref{l:propertiesF}, a revelant question is whether nonlinear resistive circuits can approximate any continuous function \emph{that satisfies Properties~\ref{i:xones} and \ref{i:xy}}. It turns out that the answer is affirmative, as made precise in the following theorem. 

\begin{theorem}
\label{t:approximation}
For any $\epsilon > 0$, compact set $\calC \subseteq \R^{N_I}$ and continuous function $f : \R^{N_I} \to \R^{N_O}$ satisfying 
\begin{equation*}
\begin{aligned} 
&f(x+\alpha \ones) = f(x) + \alpha \ones \qquad\quad \forall x \in \R^{N_I} \text{ and } \alpha \in \R; \\
&f(x) \preceq f(y) \qquad\qquad\qquad\quad\: \forall x,y \in \R^{N_I} \text{ s.t. } x \preceq y,
\end{aligned}
\end{equation*}
there exists a function $F \in \calF_{N_I,N_O}$ satisfying
\begin{equation*}
\|f(x) - F(x)\| \leq \epsilon 
\end{equation*}
for all $x \in \calC$. 
\end{theorem}

The remainder of this section is devoted to proving Theorem~\ref{t:approximation}. The proof is technical, and will be split up in several parts. First, in Section~\ref{sec:singleoutput}, we will show how the problem can be reduced to the case with a single output, allowing us to focus on the proof of Theorem~\ref{t:approximation} for $N_O = 1$ in the remainder. Thereafter, in Section~\ref{sec:minmax} we will present a min-max approximation of $f$, inspired by results on min-max neural networks \cite{Sill1997,Minin2010} approximating monotone functions. Such a min-max approximation, however, cannot be directly implemented by a resistive circuit. Therefore, in Section~\ref{sec:circuitapprox}, we will show how to construct a circuit that implements a smooth approximation of the min-max function. As we will see, this requires a scaling of some of the circuit's co-content functions. 

\subsection{Reduction to the single-output case}
\label{sec:singleoutput}

The purpose of this subsection is to show that it suffices to prove Theorem~\ref{t:approximation} for the case $N_O = 1$. We will do this by proving, in Lemma~\ref{l:twographs}, that the output of a circuit composed of two subgraphs with the same set of input nodes is related to the outputs of the subcircuits. This will essentially allow us to consider the function approximation of each entry $f_i(x)$ by a subcircuit, for every $i = 1,2,\dots,N_O$. 

To facilitate the discussion, let $\calG_1 = (\calV_1,\calE_1)$ and $\calG_2 = (\calV_2,\calE_2)$ be two connected undirected graphs with the same set of input nodes, i.e., 
\begin{equation*}
\calV_i = \calV_I \cup \calV_{H_i} \cup \calV_{O_i},
\end{equation*}
where $\calV_I$, $\calV_{H_i}$ and $\calV_{O_i}$ are mutually disjoint for $i = 1,2$. Label the edges of $\calG_1$ by $1,2,\dots,|\calE_1|$ and those of $\calG_2$ by $|\calE_1|+1,|\calE_1|+2,\dots,B$, where $B:=|\calE_1|+|\calE_2|$. Let $G_k$ be co-content functions satisfying the assumption in \ref{i:func} for all $k = 1,2,\dots,B$. Let $F^1(\cdot)$ be the input-ouput map of the circuit defined by the graph $\calG_1$ and the co-content functions $G_k$ for $k = 1,2,\dots,|\calE_1|$. Let $F^2(\cdot)$ be the input-output map of the circuit defined by the graph $\calG_2$ and the functions $G_k$ for $k=|\calE_1|+1,|\calE_1|+2,\dots,B$.

Now, define the new graph $\calG = (\calV,\calE)$ by 
\begin{equation*}
\calV := \calV_I \cup \calV_{H_1} \cup \calV_{H_2} \cup\calV_{O_1} \cup \calV_{O_2}, \text{ and } \calE = \calE_1 \cup \calE_2,
\end{equation*}
and the total co-content
\begin{equation*}
G(v) := \sum_{k = 1}^{B} G_k(v_k).
\end{equation*}
Finally, denote by $F(\cdot)$ the input-output map of the circuit corresponding to the graph $\calG$ and functions $G_k$ for $k = 1,2,\dots,B$. This leads to the following lemma. 

\begin{lemma}
\label{l:twographs}
With $F,F^1$ and $F^2$ as above, we have that
\begin{equation}
\label{eqtwographs}
F(x) = \begin{bmatrix}
F^1(x) \\ F^2(x) 
\end{bmatrix}
\end{equation} 
for all $x \in \R^{|\calV_I|}$.
\end{lemma}

\begin{proof}
The incidence matrix $D$ of $\calG$ takes the form
\begin{equation}
\label{Dstructured}
D = \begin{pmat}[{.}]
D_I^1 & D_I^2 \cr\- D_H^1 & 0 \cr 0 & D_H^2 \cr\- D_O^1 & 0 \cr 0 & D_O^2 \cr
\end{pmat},
\end{equation}
where 
\begin{equation*}
\begin{bmatrix}
D_I^i \\ D_H^i \\ D_O^i
\end{bmatrix}
\end{equation*}
is the incidence matrix of $\calG_i$ for $i = 1,2$. Let $x \in \R^{|\calV_I|}$ and partition 
\begin{equation*}
F(x) = \begin{bmatrix}
F_1(x) \\ F_2(x) 
\end{bmatrix},
\end{equation*}
where $F_i(x) \in \R^{|\calV_{O_i}|}$ for $i = 1,2$. Let $p_{H_1}$ and $p_{H_2}$ be such that
\begin{equation*}
\begin{bmatrix}
D_H^1 & 0 \\ 0 & D_H^2 \\ D_O^1 & 0 \\ 0 & D_O^2
\end{bmatrix} \nabla G\left( D^\top \begin{bmatrix}
x \\ p_{H_1} \\ p_{H_2} \\ F_1(x) \\ F_2(x)
\end{bmatrix} \right) = 0.
\end{equation*}
Due to the fact that
\begin{equation*}
\nabla G(v) = \begin{bmatrix}
\frac{dG_1}{dv_1}(v_1) \\ \frac{dG_2}{dv_2}(v_2) \\ \vdots \\ \frac{dG_{B}}{dv_{B}}(v_{B})
\end{bmatrix},
\end{equation*}
and the structure of the incidence matrix in \eqref{Dstructured}, we obtain, for $i = 1,2$,
\begin{equation}
\label{Fix}
\begin{bmatrix}
D_H^i \\ D_O^i
\end{bmatrix} \nabla G^i((D_I^i)^\top x + (D_H^i)^\top p_{H_i} + (D_{O}^i)^\top F_i(x)) = 0,
\end{equation} 
where
\begin{equation*}
G^1(v^1) := \sum_{k=1}^{|\calE_1|} G_k(v_k), \quad G^2(v^2) := \sum_{k=|\calE_1| + 1}^{B} G_k(v_k),
\end{equation*}
while 
\begin{equation*}
v^1 := \begin{bmatrix}
v_1 \\ v_2 \\ \vdots \\ v_{|\calE_1|}
\end{bmatrix}, \:\: \text{ and } \:\: v^2 := \begin{bmatrix}
v_{|\calE_1|+1} \\ v_{|\calE_1|+2} \\ \vdots \\ v_{B}
\end{bmatrix}.
\end{equation*}
Finally, note that \eqref{Fix} simply corresponds to the equation of the circuit defined on the graph $\calG_i$. Due to the uniqueness of the output potentials, we see that $F_i(x) = F^i(x)$ for $i = 1,2$, and thus, \eqref{eqtwographs} holds. 
\end{proof}

Lemma~\ref{l:twographs} shows that it suffices to prove Theorem~\ref{t:approximation} for $N_O = 1$. Indeed, suppose that Theorem~\ref{t:approximation} holds for $N_O = 1$. Now, let $f: \R^{N_I} \to \R^{N_O}$ satisfy the assumptions of Theorem~\ref{t:approximation} for some $N_O \in \N$. Then, for all $i = 1,2,\dots,N_O$, $f_i$ is continuous, $f_i(x + \alpha \ones) = f_i(x) + \alpha$, and $f_i(x) \leq f_i(y)$ for all $x\in \R^{N_I}$, $\alpha \in \R$ and $y \in \R^{N_I}$ satisfying $x \preceq y$. Let $\epsilon > 0 $ and let $\calC \subseteq \R^{N_I}$ be compact. Then, applying Theorem~\ref{t:approximation} with $N_O = 1$, there exist resistive circuits with input-output maps $F^i \in \calF_{N_I,1}$ such that
\begin{equation}
\label{ineqentryi}
| f_i(x) - F^i(x) | \leq \frac{\epsilon}{\sqrt{N_O}}
\end{equation} 
for all $i = 1,2,\dots,N_O$. Repeated application of Lemma~\ref{l:twographs} shows the existence of a resistive circuit with input-output map $F \in \calF_{N_I,N_O}$ satisfying
\begin{equation*}
F(x) = \begin{bmatrix}
F^1(x) \\ F^2(x) \\ \vdots \\ F^{N_O}(x) 
\end{bmatrix}
\end{equation*}
for all $ x \in \R^{N_I}$. By \eqref{ineqentryi}, this shows that
\begin{equation*}
\| f(x) - F(x) \| \leq \sqrt{N_O \frac{\epsilon^2}{N_O}} = \epsilon,
\end{equation*}
proving Theorem~\ref{t:approximation} for the general case $N_O \in \N$. Therefore, it remains to prove Theorem~\ref{t:approximation} for the case $N_O = 1$. 

\subsection{Min-max representation}
\label{sec:minmax}

The following lemma shows that all continuous functions satisfying the physical constraints can be approximated by certain min-max functions. For a positive integer $n$, we will use the shorthand notation $[n] := \{1,\dots,n\}$.

\begin{lemma}
\label{l:minmax}
Let $\epsilon > 0$ and let $\calC \subseteq \R^{N_I}$ be compact. Consider a continuous function $f:\R^{N_I} \to \R$ satisfying 
\begin{align} 
f(x+\alpha \ones) &= f(x) + \alpha \label{lminmax1} \\
f(x) &\leq f(y),  \label{lminmax2}
\end{align}
for all $\alpha \in \R$ and $x,y \in \R^{N_I}$ such that $x \preceq y$. There exists a finite subset $\calA \subseteq \calC$ such that the function $f_\calA : \R^{N_I} \to \R$, defined by 
\begin{equation*}
f_\calA(x) := \min_{a \in \calA} \left(f(a) + \max_{i \in [N_I]} (x_i - a_i) \right),
\end{equation*}
satisfies
\begin{equation*}
| f(x) - f_\calA(x) | \leq \epsilon
\end{equation*}
for all $x \in \calC$. 
\end{lemma} 
 
\begin{proof}
Since $\calC$ is compact, there exists a finite subset $\calA \subseteq \calC$ such that, for all $x \in \calC$, 
\begin{equation}
\label{xa}
\|x-a\|_\infty \leq \frac{\epsilon}{2}
\end{equation} 
for some $a \in \calA$. 

Now, let $x \in \calC$. Then 
\begin{equation*}
x \preceq a + \left(\max_{i \in [N_I]} (x_i - a_i)\right) \ones
\end{equation*}
for all $a \in \calA$, and thus, 
\begin{equation*}
f(x) \leq f\left(a + \left(\max_{i \in [N_I]} (x_i - a_i)\right) \ones\right)
\end{equation*}
for all $a \in \calA$, due to \eqref{lminmax2}. Using \eqref{lminmax1}, we see that
\begin{equation*}
f(x) \leq f(a) + \max_{i \in [N_I]} (x_i - a_i).
\end{equation*}
Since this inequality holds for all $a \in \calA$, we conclude that $f(x) \leq f_\calA(x)$.

Next, fix an $a \in \calA$ satisfying \eqref{xa}. Then,
\begin{equation*}
a-\frac{\epsilon}{2} \ones \preceq x 
\end{equation*}
and thus, 
\begin{equation}
\label{fafx}
f(a) - \frac{\epsilon}{2} \leq f(x).
\end{equation}
Due to the definition of $f_\calA$, we have that
\begin{equation*}
\begin{aligned} 
f_\calA(x) &\leq f(a) + \max_{i \in [N_I]}(x_i-a_i) \\
&\leq f(a) + \frac{\epsilon}{2} \\
&\leq f(x) + \epsilon,
\end{aligned}
\end{equation*} 
where the last inequality follows from \eqref{fafx}. We conclude that $0 \leq f_\calA(x) - f(x) \leq \epsilon$, and therefore, $| f_\calA(x) - f(x) | \leq \epsilon$. This proves the lemma. 
\end{proof}

\subsection{Circuit approximation of min-max and proof of Theorem~\ref{t:approximation}}
\label{sec:circuitapprox}

In view of Lemma~\ref{l:minmax}, the idea is now to build a circuit that approximates the function $f_\calA(x)$. To do so, we will make use of smooth approximations of the $\max$ and $\min$ functions. For $\lambda > 0$ and $x \in \R^n$, we define
\begin{equation*}
\begin{aligned} 
\smax_\lambda(x) &:= \frac{1}{\lambda} \ln \left( \frac{1}{n} \sum_{i=1}^n e^{\lambda x_i} \right), \text{ and } \\
\smin_\lambda(x) &:= -\frac{1}{\lambda} \ln \left( \frac{1}{n} \sum_{i=1}^n e^{-\lambda x_i} \right).
\end{aligned}
\end{equation*} 

These functions are scaled versions of the LogSumExp function, commonly encountered in machine learning applications, see especially the relevant work \cite{Igel2023}. The following basic lemma clarifies in which sense these functions are approximations of $\max$ and $\min$. The lemma also asserts that a property of the $\max$ and $\min$ functions is carried over to $\smax_\lambda$ and $\smin_\lambda$, respectively. More specifically, the $\max$ and $\min$ functions can be shown to be nonexpansive with respect to the $\infty$-norm, which also holds true for $\smax_\lambda$ and $\smin_\lambda$.

\begin{lemma}
\label{l:sminmax}
For any $\lambda > 0$ and $x \in \R^n$, we have that
\begin{equation*}
\begin{aligned} 
|\smax_\lambda(x) - \max_{i \in [n]} x_i | &\leq \frac{\ln(n)}{\lambda}, \text{ and } \\
|\smin_\lambda(x) - \min_{i \in [n]} x_i | &\leq \frac{\ln(n)}{\lambda}.
\end{aligned}
\end{equation*} 

\noindent
Furthermore, both $\smax_\lambda$ and $\smin_\lambda$ are nonexpansive with respect to the $\infty$-norm, that is, 
\begin{align}
|\smax_\lambda(x) - \smax_\lambda(y)| &\leq \|x-y\|_\infty \label{smaxLC}\\
|\smin_\lambda(x) - \smin_\lambda(y)| &\leq \|x-y\|_\infty, \label{sminLC}
\end{align}
for all $x,y \in \R^n$.
\end{lemma} 
The proof of Lemma~\ref{l:sminmax} is given in Appendix~\ref{sec:prooflemsminmax}.

We are now in a position to prove Theorem~\ref{t:approximation}. Roughly speaking, the idea is as follows: we will construct a circuit with one layer of input nodes, one layer of hidden nodes and one output node. The purpose of the potentials of the hidden nodes is to approximate the ``inner" maximum defining $f_\calA$, while the output potential approximates the ``outer" minimum in the definition of $f_\calA$. A technical challenge, addressed in the proof of Theorem~\ref{t:approximation}, is that the hidden potentials cannot compute the smooth approximation $\smax_\lambda$ of the max-function exactly, due to the presence of branches between hidden nodes and the output node. As we will see, this is overcome by a suitable scaling of some of the co-content functions in the network.

\begin{proof}[Proof of Theorem~\ref{t:approximation}]
As shown in Section~\ref{sec:singleoutput}, it suffices to prove the result for $N_O = 1$. Let $\epsilon > 0$, $\calC \subseteq \R^{N_I}$ be compact, and let $f : \R^{N_I} \to \R$ be a function satisfying the assumptions of Theorem~\ref{t:approximation}. By Lemma~\ref{l:minmax}, there exists a finite subset $\calA \subseteq \calC$ such that 
\begin{equation}
\label{diffxfA}
|f(x) - f_\calA(x) | \leq \frac{\epsilon}{3}
\end{equation} 
for all $x \in \calC$. Let $q:=|\calA|$ and denote by $a^1,a^2,\dots,a^q$ the elements of $\calA$. Let $\lambda > 0$, whose precise value will be determined later, and define $\bar{z} : \R^{N_I} \to \R^{q}$ elementwise by
\begin{equation}
\label{defzj}
\bar{z}_j(x) := \smax_\lambda(f(a^j) \ones + x - a^j),
\end{equation} 
for $x \in \calC$ and $j \in [q]$. Also define
\begin{equation*}
\bar{y}(x) := \smin_\lambda(\bar{z}(x))
\end{equation*}
for $x \in \R^{N_I}$. Now, by the triangle inequality,
\begin{equation*}
\begin{aligned}
&|\bar{y}(x) - f_\calA(x)| \leq
|\smin_\lambda(\bar{z}(x)) - \min_{j \in [q]} \bar{z}_j(x)| \: + \\
&|\min_{j \in [q]} \bar{z}_j(x) - \min_{j \in [q]} \max_{i \in [N_I]} (f(a^j) + x_i - a_i^j)|.
\end{aligned}
\end{equation*}
The first term on the right hand side is upper bounded by $\frac{\ln(q)}{\lambda}$, due to Lemma~\ref{l:sminmax}. By nonexpansiveness of $\min$ with respect to the $\infty$-norm, the second term is upper bounded by 
$$\left\| \bar{z}(x) - \begin{bmatrix}
    \max_{i \in [N_i]} (f(a^1) + x_i-a_i^1) \\ \vdots \\ \max_{i \in [N_i]} (f(a^q)+ x_i-a_i^q) \end{bmatrix} \right\|_\infty,
$$
which, by definition of the $\infty$-norm, equals
$$
|\bar{z}_k(x) - \max_{i \in [N_I]} (f(a^k) + x_i - a_i^k)|
$$
for some index $k \in [q]$. By Lemma~\ref{l:sminmax}, the latter term is upper bounded by $\frac{\ln(N_I)}{\lambda}$. We conclude that 
$$
|\bar{y}(x) - f_\calA(x)| \leq \frac{\ln(q)}{\lambda} + \frac{\ln(N_I)}{\lambda}
$$
for all $x \in \calC$.

Now, choose $\lambda > 0$ sufficiently large such that
\begin{equation}
\label{difbaryfA}
|\bar{y}(x) - f_\calA(x)| \leq \frac{\epsilon}{3}
\end{equation}
for all $x \in \calC$. 

Next, consider a connected graph $\calG = (\calV,\calE)$, where $\calV = \calV_I \cup \calV_H \cup \calV_O$ with $|\calV_I| = N_I$, $|\calV_H| = q$ and $|\calV_O| = 1$. The set $\calE$ consists of $(N_I+1)q$ branches: each input node in $\calV_I$ is connected to all hidden nodes in $\calV_H$, and all hidden nodes are connected to the output node. There are no branches between input nodes and the output node.

We will now define the co-content functions associated to the branches. For each input node $i = 1,2,\dots,N_I$ and hidden node $j = 1,2,\dots,q$, we consider, for $k = (i-1)q+j$, the function
\begin{equation}
\label{cocont1}
G_k(v_k) := c \left( v_k + \frac{1}{\lambda} e^{-\lambda (v_k - f(a^j) + a_i^j)} \right)
\end{equation}
as the co-content of the branch between $i$ and $j$. 
Here $c > 0$ is a to-be-determined constant. Furthermore, define, for $k = N_I q + j$ with $j = 1,2,\dots,q$, the function 
\begin{equation}
\label{cocont2}
G_k(v_k) := \frac{1}{\lambda} e^{\lambda v_k} - v_k
\end{equation}
as the co-content of the branch between the hidden node $j$ and the output node. Note that the functions $G_k$ indeed satisfy the assumptions in \ref{i:func} for all $k = 1,2,\dots,(N_I+1)q$. The total co-content function is now given by
\begin{equation*}
G(v) = \sum_{k = 1}^{(N_I + 1)q} G_k(v_k).
\end{equation*}
After assigning an orientation to each branch, the incidence matrix of the network takes the form
\begin{equation}
\label{incuni}
D = \begin{bmatrix}
D_I \\ D_H \\ D_O
\end{bmatrix} = \begin{bmatrix}
-e_1 \ones^\top & \cdots & -e_{N_I} \ones^\top & 0 \\
I & \cdots & I & -I \\
0 & \cdots & 0 & \ones^\top
\end{bmatrix},
\end{equation}
where $e_i$ denotes the $i$-th standard basis vector of $\R^{N_I}$ for $i = 1,2,\dots,N_I$. We then have that 
\begin{equation*}
\begin{aligned}
&G(D_I^\top x + D_H^\top z + D_O^\top y) = \sum_{j=1}^q \frac{1}{\lambda} e^{\lambda (y-z_j)}-(y-z_j) \\ &+c \sum_{j = 1}^q \sum_{i = 1}^{N_I} \left( z_j - x_i + \frac{1}{\lambda} e^{-\lambda(z_j-x_i-f(a^j)+a_i^j)} \right),
\end{aligned}
\end{equation*}
for all $x \in \R^{N_I}$, $z \in \R^q$ and $y \in \R$. Now, define the functions
\begin{equation*}
\begin{aligned}
H_1(x,z) &:= \sum_{j = 1}^q \sum_{i = 1}^{N_I} \left( z_j - x_i + \frac{1}{\lambda} e^{-\lambda(z_j-x_i-f(a^j)+a_i^j)} \right) \\ 
H_2(z,y) &:= \sum_{j=1}^q \frac{1}{\lambda} e^{\lambda (y-z_j)}-(y-z_j),
\end{aligned}
\end{equation*}
so that 
\begin{equation*}
G(D_I^\top x + D_H^\top z + D_O^\top y) = cH_1(x,z) + H_2(z,y).
\end{equation*}
Given $c > 0$ and the input potentials $x \in \R^{N_I}$, denote by $(z^c(x),y^c(x))$ the corresponding hidden and output potentials, i.e., 
\begin{equation*}
(z^c(x),y^c(x)) = \argmin_{(z,y) \in \R^{q} \times \R} G(D_I^\top x + D_H^\top z + D_O^\top y).
\end{equation*}
Existence and uniqueness of these potentials follows from Lemma~\ref{l:uniqueoutput}. We emphasize that $z^c(x)$ is not equal to $\bar{z}(x)$, i.e., the hidden potentials $z^c(x)$ do not compute the $\smax_\lambda$ functions in \eqref{defzj} exactly. We will, however, show that $z^c(\cdot)$ converges uniformly to $\bar{z}(\cdot)$ on $\calC$ as $c \to \infty$. To do so, we need some preparation. 

First, note that
\begin{equation*}
\begin{aligned}
&cH_1(x,z^c(x)) + H_2(z^c(x),y^c(x)) \leq \\ &c H_1(x,\bar{z}(x)) + H_2(\bar{z}(x),y),
\end{aligned}
\end{equation*}
for all $x \in \R^{N_I}$ and $y \in \R$. It can be readily shown that the minimum of the function $H_2$ is $\frac{q}{\lambda}$, and thus,
\begin{equation*}
H_2(z^c(x),y^c(x)) \geq \frac{q}{\lambda}.
\end{equation*}
Furthermore, the function $x \mapsto H_2(\bar{z}(x),y)$ is continuous for all $y$. Fix $y \in \R$. Since $\calC$ is compact and $H_2$ and $\bar{z}$ are continuous, there exists a constant $a > 0$ such that
\begin{equation*}
H_2(\bar{z}(x),y) \leq a
\end{equation*}
for all $x \in \calC$. Therefore,
\begin{equation*}
\begin{aligned}
c(H_1(x,z^c(x)) - H_1(x,\bar{z}(x))) &\leq \\  H_2(\bar{z}(x),y) - H_2(z^c(x),y^c(x)) &\leq a - \frac{q}{\lambda} =: \varphi
\end{aligned}
\end{equation*}
for all $x \in \calC$. We conclude that
\begin{equation*}
H_1(x,z^c(x)) - H_1(x,\bar{z}(x)) \leq \frac{\varphi}{c}
\end{equation*}
for all $x \in \calC$. We claim that, in addition, $H_1(x,z^c(x)) - H_1(x,\bar{z}(x)) \geq 0$. We will prove this by showing that $\bar{z}(x)$ is the minimizer of $H_1(x, \cdot)$ for all $x \in \calC$. First, note that $H_1(x, \cdot)$ is continuous, coercive and strictly convex and thus has a unique minimizer for all $x \in \calC$ by Proposition~\ref{prop:basicconvexity} in Appendix~\ref{app:convexity}. To find this minimizer, fix $x \in \calC$, and compute
\begin{equation*}
\left(\frac{\partial H_1}{\partial z}(x,z)\right)_j = \sum_{i=1}^{N_I} (1 - e^{-\lambda(z_j - x_i -f(a^j) + a_i^j)}).
\end{equation*}
Setting this equal to zero yields, for $j= 1,2,\dots,q$, 
\begin{equation*}
e^{-\lambda z_j} \sum_{i=1}^{N_I} e^{\lambda(f(a^j) + x_i - a_i^j)} = N_I.
\end{equation*}
In other words, 
\begin{equation*}
e^{\lambda z_j} = \frac{1}{N_I} \sum_{i=1}^{N_I} e^{\lambda(f(a^j) + x_i - a_i^j)},
\end{equation*}
hence, by taking natural logarithms,
\begin{equation*}
\begin{aligned}
z_j &= \frac{1}{\lambda} \ln\left( \frac{1}{N_I} \sum_{i=1}^{N_I} e^{\lambda(f(a^j) + x_i - a_i^j)} \right) \\
&= \smax_\lambda(f(a^j)\ones + x - a^j) = \bar{z}_j(x).
\end{aligned}
\end{equation*} 
So indeed, $\bar{z}(x)$ is the unique minimizer of $H_1(x, \cdot)$. We conclude that our claim holds, and thus, 
\begin{equation}
\label{uplowH1}
0 \leq H_1(x,z^c(x)) - H_1(x,\bar{z}(x)) \leq \frac{\varphi}{c}
\end{equation}
for all $x \in \calC$. 

We will now show that the function $z^c(\cdot)$ converges uniformly to $\bar{z}(\cdot)$ on $\calC$, as $c \to \infty$. To this end, fix $\rho > 0$. Define $\alpha(\rho)$ as
\begin{equation*}
\inf_{x \in \calC, \: \zeta \in \R^q} \left\{ H_1(x,\zeta) - H_1(x,\bar{z}(x)) \mid \|\zeta-\bar{z}(x)\| \geq \rho \right\}.
\end{equation*}
Due to the compactness of $\calC$ and the coercivity of $H_1$, for every $x \in \calC$, $H_1(x,\zeta) - H_1(x,\bar{z}(x)) \to \infty$ as $\| \zeta \| \to \infty$. Hence, there exists a sufficiently large constant $R > 0$ such that $\alpha(\rho)$ equals
\begin{equation*}
\inf_{x \in \calC, \: \|\zeta \| \leq R} \left\{ H_1(x,\zeta) - H_1(x,\bar{z}(x)) \mid \|\zeta-\bar{z}(x)\| \geq \rho \right\}.
\end{equation*}
Now, since $\bar{z}(x)$ is the unique minimizer of $H_1(x, \cdot)$ for all $x \in \calC$, we have that 
\begin{equation*}
H_1(x,\zeta) - H_1(x,\bar{z}(x)) > 0
\end{equation*}
for all $\zeta \in \R^q$ satisfying $ \|\zeta-\bar{z}(x)\| \geq \rho > 0$. Therefore, since
\begin{equation*}
\{(x,\zeta) \in \calC \times \R^q \mid \| \zeta \| \leq R \}
\end{equation*}
is compact and $H_1$ and $\bar{z}$ are continuous, we also have that $\alpha(\rho) > 0$. 

Next, choose $c > 0$ sufficiently large so that 
\begin{equation*}
\frac{\varphi}{c} < \alpha(\rho).
\end{equation*}
It then follows from \eqref{uplowH1} that 
\begin{equation*}
H_1(x,z^c(x)) - H_1(x,\bar{z}(x)) < \alpha(\rho)
\end{equation*}
for all $x \in \calC$. We conclude that $\|z^c(x) - \bar{z}(x) \| \!<\! \rho$ for all $x \in \calC$. This proves that $z^c(\cdot)$ converges uniformly to $\bar{z}(\cdot)$ on $\calC$ as $c \to \infty$. 

The final step of the proof is to show that also $y^c(\cdot)$ converges uniformly to $\bar{y}(\cdot)$ on $\calC$ as $c \to \infty$. To this end, by the network model, we have that
\begin{equation*}
\begin{bmatrix}
D_H \\ D_O
\end{bmatrix} \nabla G(D_I^\top x + D_H^\top z^c(x) + D_O^\top y^c(x)) = 0.
\end{equation*}
By focusing on the last row, and using the structure of $D_O$ in \eqref{incuni} and the definitions of $G_k$ in \eqref{cocont2}, we see that
\begin{equation*}
\sum_{j=1}^q (e^{\lambda(y^c(x)-z_j^c(x))} -1) = 0,
\end{equation*}
and thus,
\begin{equation*}
e^{\lambda y^c(x)} \sum_{j=1}^q e^{-\lambda z_j^c(x)} = q.
\end{equation*}
This implies that
\begin{equation*}
e^{-\lambda y^c(x)} = \frac{1}{q} \sum_{j=1}^q e^{-\lambda z_j^c(x)},
\end{equation*}
showing that, after taking natural logarithms on both sides, 
\begin{equation*}
y^c(x) = - \frac{1}{\lambda} \ln \left( \frac{1}{q} \sum_{j=1}^q e^{-\lambda z_j^c(x)} \right) = \smin_\lambda(z^c(x)).
\end{equation*}
Now, recall that $\bar{y}(x) = \smin_\lambda(\bar{z}(x))$, by definition. By nonexpansiveness of $\smin_\lambda$ (see Lemma~\ref{l:sminmax}),
\begin{equation*}
\begin{aligned}
|y^c(x) - \bar{y}(x)| &= | \smin_\lambda(z^c(x)) - \smin_\lambda(\bar{z}(x))| \\ &\leq \|z^c(x) - \bar{z}(x)\|_\infty
\end{aligned}
\end{equation*}
for all $x \in \calC$. This shows that $y^c(\cdot)$ converges uniformly to $\bar{y}(\cdot)$ on $\calC$ as $c \to \infty$. 

Finally, we assemble all ingredients and complete the proof. Choose $c> 0$ sufficiently large so that
\begin{equation}
\label{difycbary}
|y^c(x) - \bar{y}(x)| \leq \frac{\epsilon}{3}
\end{equation} 
for all $x \in \calC$. Then, by combining \eqref{difycbary}, \eqref{difbaryfA}, and \eqref{diffxfA}, for all $x \in \calC$, 
\begin{equation*}
\begin{aligned}
|y^c(x) - f(x)| &\leq \\ |y^c(x) - \bar{y}(x)| + |\bar{y}(x) - f_\calA(x)| + |f_\calA(x) - f(x)| &\leq \\ \frac{\epsilon}{3} + \frac{\epsilon}{3} + \frac{\epsilon}{3},
\end{aligned}
\end{equation*}
which is equal to $\epsilon$. Thus, taking $F = y^c$ proves the theorem. 
\end{proof}

\vspace{-15pt}

\section{Conclusions and discussion}
\label{sec:conclusions}

In this paper, we have established a rigorous mathematical foundation for both the learning dynamics and the representational capacity of nonlinear resistive networks with monotone elements. By framing the solution of these circuits as the global minimizer of a strictly convex co-content function, we proved that Target-Clamped Contrastive Learning (TCCL) reliably converges to parameters that implement desired input-output mappings. 

Crucial to our analysis has been a detailed study of the properties of the contrastive function, as reported in Thms.~\ref{prop:Q} and \ref{t:convexQ}. 

In the case of only input and output nodes, the contrastive function is convex (Thm.~\ref{t:convexQ}), even for networks with non-monotone resistors (Rem.~\ref{rem:convexnonmonotone}). As a consequence, TCCL converges to a global minimizer of the contrastive function for a range of stepsizes (Thm.~\ref{t:convergenceconvex}). This property mirrors the convexity of linear regression and logistic regression tasks in traditional machine learning. In contrast, we have shown that the objective function used in Equilibrium Propagation (EP) is in general \emph{not convex}, even for linear resistive circuits (Prop.~\ref{prop:Jnotconvex}). This demonstrates an important advantage of using TCCL in networks without hidden nodes. A further advantage of TCCL over EP and coupled learning is that it does not involve a nudging parameter, thereby eliminating the need to fine-tune it, and simplifying the learning procedure altogether. 

Beyond the setting with only input and output nodes, we have revealed general properties of the contrastive function, also in the presence of hidden nodes (Thm.~\ref{prop:Q}). In this setting, the contrastive function is generally not convex, but it exhibits an attractive scaling behavior. Assuming all resistors are trainable, the hidden and output potentials of a resistive network are invariant under scaling of the parameters (Lem.~\ref{l:scalingpO}), and only depend on the relative magnitudes of parameters. As noted in \cite{Stern2024powerefficient}, this invariance implies that network parameters can be uniformly rescaled to reduce operating power without altering the input-output mapping. As we have shown here, this feature implies that the contrastive function is positively homogeneous of degree~$1$ (Thm.~\ref{prop:Q}). Its minimization therefore has a regularization effect: it penalizes parameter configurations with large parameter values, and thereby promotes low (pseudo-)power dissipation during inference. This is achieved without introducing an additional term in the learning rule, as is done in \cite{Stern2024powerefficient}. In contrast, the EP objective is positively homogeneous of degree~$0$ (Lem.~\ref{l:scalingpO}), and consequently exhibits no inherent preference between high- and low-power parameter configurations. In the presence of hidden nodes, we have shown that TCCL converges to a stationary point on a subsequence (Thm.~\ref{t:convergencehidden}). This is in line with similar results on the convergence of (stochastic) gradient descent applied to deep neural networks, see e.g., \cite[Sec. 4.3]{Bottou2018}. However, different from deep neural networks, due to the positive homogeneity of the contrastive function, all limit points in the interior of the constraint set are \emph{global minimizers} (Thm.~\ref{t:convergencehidden}).

Furthermore, we characterized the computational expressivity of these resistive networks. We demonstrated that their underlying physics inherently restricts them to monotonically non-decreasing and shift-equivariant functions (Lem.~\ref{l:propertiesF}). Crucially, we established a novel approximation theorem proving that these constraints are the only limitations on expressivity, as any continuous function satisfying these two properties can be approximated to arbitrary precision by a sufficiently large network with monotone elements (Thm.~\ref{t:approximation}).

Beyond the electrical domain, the framework developed in this work has broader implications for other analog learning systems. Our results generalize naturally to any physical network governed by analogous flow-effort relationships. For instance, mechanical elastic networks characterized by nonlinear force-displacement relationships minimize potential energy in a manner that is similar to co-content minimization in resistive circuits \cite[Sec. II.B]{Scellier2026quantum}. Consequently, the convergence guarantees and approximation results established here can be directly translated to design, train, and analyze such networks. Our learning framework may also be considered for coupled phase oscillators, for which EP has been employed in \cite{Wang2024}. In case the drive phases are fixed (non-trainable), the energy function of \cite[Eq. (1)]{Wang2024} is linear in the parameters, and thus, the contrastive function is degree-$1$ positively homogeneous. However, a challenge compared to the current work is the non-convexity of the energy as a function of the phase angles.

The primary limitation of this study is the restriction to strictly monotone elements, which limits the class of functions the network can approximate. This constraint can be overcome by introducing active two-port elements, such as VCVSs \citep{Scellier2025}, or by incorporating non-monotone components (with regions of negative differential conductance), similar to the approach of \cite{lipka2026thermodynamic}. However, by introducing non-monotone elements, the network's energy function (co-content) is no longer convex. For general non-monotone circuit elements, the energy landscape can develop \emph{multiple} stationary points of different types, each corresponding to a valid physical solution. This multiplicity renders the circuit's input-output mapping mathematically ill-defined. Developing the mathematical tools to rigorously study expressivity and learning in networks with multiple stationary points remains a key direction for future work. In these networks, EP may prove to be a preferable training strategy, as EP relies only on stationarity rather than global energy minimization; see for instance \cite{kerjan2026training,hoier2026training} where EP was successfully employed to (numerically) train other types of energy-based models to classify ImageNet and reconstruct corrupted images.

Another promising direction for future research concerns the representational capacity of networks consisting exclusively of input and output nodes. Thus far, we have investigated the approximation properties of monotone resistive circuits with hidden nodes (Thm.~\ref{t:approximation}). However, in view of the strong convergence guarantees established for networks without hidden nodes (Thm.~\ref{t:convergenceconvex}), it is relevant to characterize the functions that such networks can implement, and to compare these to the class of linear and logistic classifiers in traditional machine learning. A further direction is to extend the framework beyond static resistive circuits by incorporating dynamic circuit elements, such as capacitors, thereby enabling learning of dynamical systems. An initial step in this direction has already been taken in the recent work \cite{Heidema2026}, which extends target-clamped contrastive learning to linear resistive-capacitive networks.

In this work, we have assumed that the data consist of a single pair of input-output potentials. This assumption was made for simplicity and clarity of presentation, and we expect that the framework can be extended straightforwardly to accommodate multiple data points. Indeed, such an extension has already been established in \cite[Sec. V]{huijzer2026convergence} for linear networks without hidden nodes. Here, the idea was to extend the convergence analysis of projected gradient descent for a single data sample to that of \emph{stochastic} projected gradient descent, thereby allowing for multiple samples. We also note that, while our theoretical analysis is restricted to optimization and expressivity, investigating the generalization performance of nonlinear resistive networks on unseen data is an important avenue for future research.

Finally, we close this paper with a few remarks regarding the origins of our work and the employed terminology. Contrastive Learning (CL) finds its roots in the study of Boltzmann machines \cite{ackley1985learning}, where target clamping was shown to minimize the Kullback-Leibler divergence between the data and model distributions. CL with target clamping was later extended to continuous Hopfield networks \cite{Movellan1991}, but in recent years, the physical learning community has favored nudging-based approaches such as EP and coupled learning, while keeping ``contrastive learning” as an umbrella term for the broader family of energy-based learning methods \cite{Stern2023annualreview,momeni2025training}. Given the theoretical advantages of fully clamping outputs to their target values, as demonstrated in this work and other studies \cite{huijzer2026convergence,Heidema2026}, we believe this approach deserves further research. We therefore adopt the term Target-Clamped Contrastive Learning (TCCL) to distinguish it from nudging-based approaches. TCCL may also be studied in settings where EP has previously been adapted, including Lagrangian dynamical systems \cite[Sec. 6.1]{scellier2021deep} and quantum systems \cite{Scellier2026quantum}, where the principle of co-content minimization is replaced by the principle of least action and the variational principle of quantum mechanics, respectively. Similarly, TCCL may be transposed to bilevel optimization problems \cite{zucchet2022beyond}.

\appendix

\section{Notation regarding partial derivatives}
\label{app:notation}

We first introduce some notation related to partial derivatives. For a differentiable function $f: \R^n \to \R$ we use the notation
\begin{equation*}
\frac{\partial f}{\partial x}(x) := \begin{bmatrix}
\frac{\partial f}{\partial x_1}(x) \\ 
\frac{\partial f}{\partial x_2}(x) \\ 
\vdots \\
\frac{\partial f}{\partial x_n}(x) 
\end{bmatrix}, 
\end{equation*}
and 
\begin{equation*}
\frac{\partial f}{\partial x^\top}(x) := \begin{bmatrix}
\frac{\partial f}{\partial x_1}(x) & 
\frac{\partial f}{\partial x_2}(x) &
\cdots &
\frac{\partial f}{\partial x_n}(x) 
\end{bmatrix}.
\end{equation*}
We will also use the shorthand notation $\nabla f := \frac{\partial f}{\partial x}$ to denote the gradient of $f$. The Jacobian of a differentiable function $g: \R^n \to \R^m$ is defined as
\begin{equation*}
\frac{\partial g}{\partial x^\top}(x) := \begin{bmatrix}
\frac{\partial g_1}{\partial x_1}(x) & \cdots & \frac{\partial g_1}{\partial x_n}(x) \\ 
\vdots & \ddots & \vdots \\
\frac{\partial g_m}{\partial x_1}(x) & \cdots & \frac{\partial g_m}{\partial x_n}(x)
\end{bmatrix}.
\end{equation*}
The Hessian of a twice differentiable $f: \R^n \to \R$ is:
\begin{equation*}
\begin{aligned} 
&\frac{\partial^2 f}{\partial x \partial x^\top}(x) := \frac{\partial}{\partial x^\top}\left(\frac{\partial f}{\partial x}\right)(x) \\ &= \begin{bmatrix}
\frac{\partial^2 f}{\partial x_1^2}(x) & \frac{\partial^2 f}{\partial x_1 \partial x_2}(x) & \cdots & \frac{\partial^2 f}{\partial x_1 \partial x_n}(x) \\
\frac{\partial^2 f}{\partial x_2 \partial x_1}(x) & \frac{\partial^2 f}{\partial x_2^2}(x) & \cdots & \frac{\partial^2 f}{\partial x_2 \partial x_n}(x) \\
\vdots & \vdots & \ddots & \vdots \\
\frac{\partial^2 f}{\partial x_n \partial x_1}(x) & \frac{\partial^2 f}{\partial x_n \partial x_2}(x) & \cdots & \frac{\partial^2 f}{\partial x_n^2}(x) 
\end{bmatrix}.
\end{aligned} 
\end{equation*} 
We will use the shorthand notation $\nabla^2 f := \frac{\partial^2 f}{\partial x \partial x^\top}$.

Finally, for a function $f : \R^{n_1} \times \cdots \times \R^{n_m} \to \R$ in the $m$ variables $x^1,\dots,x^m$, we use the notation
\begin{equation*} 
\begin{aligned} 
\frac{\partial f}{\partial x^i}(x^1,\dots,x^m) &:= \begin{bmatrix}
\frac{\partial f}{\partial x^i_1}(x^1,\dots,x^m) \\  \vdots \\ 
\frac{\partial f}{\partial x^i_{n_i}}(x^1,\dots,x^m)
\end{bmatrix},
\end{aligned}
\end{equation*}
$
\frac{\partial f}{\partial (x^j)^\top} := 
\left(\frac{\partial f}{\partial x^j}\right)^\top
$
and 
\begin{equation*}
\frac{\partial^2 f}{\partial x^i \partial (x^j)^\top} := 
\frac{\partial}{\partial (x^j)^\top} \left(\frac{\partial f}{\partial x_i}\right) 
\end{equation*}
for $i,j = 1,2,\dots,m$. 

\vspace{-5pt}

\section{Preliminaries on Lipschitz continuity, coercivity and convexity}
\label{app:convexity}

We call $f: D \subseteq \R^n \to \R^m$ \emph{Lipschitz continuous} if there exists a constant $L \geq 0$ such that
\begin{equation*}
\|f(x) - f(y)\| \leq L \|x - y\|
\end{equation*}
for all $x,y \in D$. We call such an $L$ a \emph{Lipschitz constant} of $f$. In case $f$ is Lipschitz continuous with $L = 1$, we call $f$ \emph{nonexpansive}. The following result for differentiable functions with Lipschitz continuous gradients is well-known (see, e.g., \cite[Lem. 5.7]{Beck2017}).
\begin{lemma}
\label{l:Lipschitzineq}
Let $f: \R^n \to \R$ be differentiable and suppose $\nabla f$ is Lipschitz continuous with constant $L$. Then
\begin{equation*}
f(x) \leq f(y) + (x-y)^\top \nabla f(y) + \frac{L}{2} \|x-y\|^2
\end{equation*}
for all $x,y \in \R^n$.
\end{lemma}

A function $f: \mathbb{R}^n \to \mathbb{R}$ is called \emph{coercive} if $f(x) \to \infty$ as $\|x\| \to \infty$.  

A subset $\mathcal{C} \subseteq \R^n$ is called \emph{convex} if for any two points $x,y \in \mathcal{C}$ we have that $\alpha x + (1-\alpha)y \in \mathcal{C}$ for all $\alpha \in [0,1]$. A function $f: \mathcal{C} \to \mathbb{R}$ is \emph{convex} if
\begin{equation*}
f(\alpha x + (1-\alpha)y) \leq \alpha f(x) + (1-\alpha) f(y)
\end{equation*}
for all $x,y \in \mathcal{C}$ and all $\alpha \in [0,1]$. We call $f$ \emph{concave} if $-f$ is convex. A convex function $f : \mathcal{C} \to \mathbb{R}$ is called \emph{strictly convex} if 
\begin{equation*}
f(\alpha x + (1-\alpha) y) < \alpha f(x) + (1-\alpha) f(y)
\end{equation*}
for all distinct $x,y \in \mathcal{C}$ and all $\alpha \in (0,1)$. If the function $f : \calC \to \R$ is twice continuously differentiable, strict convexity is implied by $\nabla^2 f(x) > 0$ for all $x \in \calC$, i.e., the Hessian is \emph{positive definite}. 

We now recall the following basic proposition about (strictly) convex and coercive functions. 

\begin{proposition}
\label{prop:basicconvexity}
Let $f : \mathbb{R}^n \to \mathbb{R}$. Then the following statements hold:
\begin{enumerate}[label=(\roman*)]
\item \label{item:convex} If $f$ is convex, then any local minimum of $f$ is also a global minimum.
\item \label{item:sconvex} If $f$ is strictly convex, then it has at most one global minimum. 
\item \label{item:coercive} If $f$ is continuous and coercive, then it has at least one global minimum. 
\item \label{item:global} Thus, if $f$ is continuous, coercive and strictly convex, it has precisely one global minimum. 
\end{enumerate}
\end{proposition}
Proposition~\ref{prop:basicconvexity}~\ref{item:convex} and \ref{item:sconvex} follow from \cite[Prop. 3.1.1]{Bertsekas2009}, while Proposition~\ref{prop:basicconvexity}~\ref{item:coercive} follows from Weierstrass' theorem, see e.g. \cite[Prop. 3.2.1]{Bertsekas2009}.

The \emph{projection} $P_\calC : \R^{n} \to \calC$ onto a nonempty, closed and convex set $\calC \subseteq \R^n$ is defined as
\begin{equation}
\label{defproj}
P_\calC(x) = \argmin_{y \in \calC} \|x-y\|.
\end{equation} 
The following lemma collects properties $P_\calC$. 

\begin{lemma}
\label{l:projection}
Let $\calC \subseteq \R^n$ be nonempty, closed and convex. Then the following properties hold:
\begin{enumerate}[label=\roman*)]
\item \label{projineq} We have that $z = P_\calC(x)$ if and only if 
\begin{equation*}
(x-z)^\top (y-z) \leq 0
\end{equation*} 
for all $y \in \calC$. 
\item \label{projfnonex} For all $x,y \in \R^n$, 
\begin{equation*}
(x-y)^\top (P_\calC(x) - P_\calC(y)) \geq \| P_\calC(x) - P_\calC(y) \|^2.
\end{equation*}
\item \label{projnonex} $P_\calC$ is nonexpansive. 
\end{enumerate}
\end{lemma}

\newpage 

Lemma~\ref{l:projection}~\ref{projineq} follows from \cite[Prop. 1.1.9]{Bertsekas2009}, while Lemma~\ref{l:projection}~\ref{projfnonex} follows from \cite[Prop. 4.8]{Bauschke2017}. This property is called \emph{firm nonexpansiveness}. Finally, Lemma~\ref{l:projection}~\ref{projnonex} simply follows from Lemma~\ref{l:projection}~\ref{projfnonex} and the triangle inequality. 

\section{Proof of Proposition~\ref{p:PGDconvex}}
\label{sec:proofpropconv}


Let $x^\ast \in \argmin_{x \in \calC} f(x)$. By \cite[Prop. 1.1.9]{Bertsekas2009}, 
\begin{equation*}
\nabla f(x^\ast)^\top (y-x^\ast) \geq 0 
\end{equation*}
for all $y \in \calC$. This implies that
\begin{equation*}
\begin{aligned}
((x^\ast-\gamma \nabla f(x^\ast))-x^\ast)^\top(y-x^\ast) &= \\
- \gamma \nabla f(x^\ast)^\top (y-x^\ast) &\leq 0
\end{aligned}
\end{equation*}
for all $y \in \calC$, showing that
\begin{equation*}
x^\ast = P_\calC(x^\ast - \gamma \nabla f(x^\ast)),
\end{equation*}
by Lemma~\ref{l:projection}~\ref{projineq}. Now, let $t \in \N$ and define
\begin{equation*}
y^t := x^t - \gamma \nabla f(x^t),
\: \text{ and } \:
y^\ast := x^\ast - \gamma \nabla f(x^\ast).
\end{equation*}
Then $x^{t+1} = P_\calC(y^t)$ and $x^\ast = P_\calC(y^\ast)$. Hence,
\begin{equation*}
\begin{aligned} 
&\| (y^t - y^\ast) - (x^{t+1} - x^\ast) \|^2 = \\ &\| y^t - y^\ast \|^2 \!-\! 2(y^t - y^\ast)^\top (x^{t+1} - x^\ast) \!+\! \| x^{t+1} - x^\ast \|^2 \!\leq \\
&\| y^t - y^\ast \|^2 - \| x^{t+1} - x^\ast \|^2,
\end{aligned}
\end{equation*} 
where the last inequality follows from firm nonexpansiveness of $P_\calC$, see Lemma~\ref{l:projection}~\ref{projfnonex}. This shows that 
\begin{equation}
\label{xt+1xast}
\| x^{t+1} - x^\ast \|^2 + \| (y^t - y^\ast) - (x^{t+1} - x^\ast) \|^2 \leq \| y^t - y^\ast \|^2.
\end{equation}
Now,
\begin{equation*}
y^t-y^\ast
=
(x^t-x^\ast)-\gamma(\nabla f(x^t)-\nabla f(x^\ast)).
\end{equation*}
Hence
\begin{equation*}
\begin{aligned}
\|y^t \! - \! y^\ast\|^2
&=\!
\|x^t -  x^\ast\|^2
 +\gamma^2\|\nabla f(x^t)-\nabla f(x^\ast)\|^2\\
&- 2\gamma (x^t-x^\ast)^\top (\nabla f(x^t)\!-\!\nabla f(x^\ast)).
\end{aligned}
\end{equation*}
By the Baillon--Haddad theorem (see, e.g. \cite[Thm. 1.1]{Bauschke2009}),
\begin{equation*}
(x^t-x^\ast)^\top (\nabla f(x^t)-\nabla f(x^\ast)) \geq 
\frac{1}{L} \|\nabla f(x^t)-\nabla f(x^\ast)\|^2,
\end{equation*}
and thus, 
\begin{equation*}
\|y^t-y^\ast\|^2 \leq \|x^t-x^\ast\|^2 - \gamma (\frac{2}{L} - \gamma) \|\nabla f(x^t)-\nabla f(x^\ast)\|^2.
\end{equation*}
Therefore, it follows from \eqref{xt+1xast} that
\begin{equation}
\label{ineqxt+1xt}
\begin{aligned}
&\| x^{t+1} - x^\ast \|^2 + \| (y^t - y^\ast) - (x^{t+1} - x^\ast) \|^2 \\ &+ \gamma (\frac{2}{L} - \gamma) \|\nabla f(x^t)-\nabla f(x^\ast)\|^2 \leq \|x^t-x^\ast\|^2.
\end{aligned}
\end{equation}
Since $\gamma \in (0,\frac{2}{L})$, this implies
\begin{equation*}
\|x^{t+1}-x^\ast\| \leq \|x^t-x^\ast\|.
\end{equation*}
Thus, $(x^t)_{t=0}^\infty$ is bounded, and, moreover, the sequence
$(\|x^t-x^\ast\|)_{t=0}^\infty$ has a limit by the monotone convergence theorem.

Summing \eqref{ineqxt+1xt} from $t = 0$ to $t = N$ yields
\begin{equation*}
\gamma (\frac{2}{L}-\gamma) \sum_{t=0}^N
\|\nabla f(x^t)-\nabla f(x^\ast)\|^2
\leq
\| x^0-x^\ast\|^2.
\end{equation*}
Hence, 
\begin{equation*}
\nabla f(x^t) - \nabla f(x^\ast) \to 0.
\end{equation*}
A similar summing argument also leads to 
\begin{equation*}
(y^t-y^\ast)-(x^{t+1}-x^\ast)\to 0.
\end{equation*}
Since
\begin{equation*}
\begin{aligned}
(y^t-y^\ast)-(x^{t+1}-x^\ast)
&= 
(x^t-x^{t+1}) \\
&- \gamma (\nabla f(x^t)-\nabla f(x^\ast)),
\end{aligned}
\end{equation*}
we obtain
\begin{equation*}
x^t-x^{t+1}\to 0.
\end{equation*}
Because $(x^t)_{t=0}^\infty$ is bounded, by Bolzano--Weierstrass it has a convergent subsequence
$(x^{t_j})_{j=0}^\infty$, where $0 \leq t_0 < t_1 < \cdots$. Denote the limit of this subsequence by $\bar{x}$. Since $\calC$ is closed, $\bar{x} \in \calC$. Also,
$x^{t_j+1}\to \bar x$, because $x^t-x^{t+1}\to 0$.

Taking the limit on both sides of
\begin{equation*}
x^{t_j+1}
=
P_\calC(x^{t_j} - \gamma \nabla f(x^{t_j})),
\end{equation*}
and using continuity of both $P_\calC$ and $\nabla f$ gives
\begin{equation*}
\bar{x} = P_\calC(\bar{x} -\gamma \nabla f(\bar{x})).
\end{equation*}
Therefore, by Lemma~\ref{l:projection}~\ref{projineq},
\begin{equation*}
\nabla f(\bar{x})^\top (y-\bar{x}) \geq 0
\end{equation*}
for all $y \in \calC$. Since $f$ is convex, it follows from \cite[Prop. 1.1.8]{Bertsekas2009} that
\begin{equation*}
\bar{x} \in \argmin_{x \in \calC}f(x).
\end{equation*}

Finally, since we have derived \eqref{ineqxt+1xt} for an arbitrary $x^\ast \in \argmin_{x \in \calC}f(x)$, we also obtain 
\begin{equation*}
\|x^{t+1}-\bar{x}\| \leq \|x^t-\bar{x}\|
\end{equation*}
for all $t \in \N$. Thus, the sequence $(a^t)_{t=0}^\infty$, defined by $a^t := \| x^t - \bar{x} \|$ for all $t \in \N$, is nonincreasing and bounded below. As such, by the monotone convergence theorem, $a^t \to \ell$ for some $\ell \geq 0$. But along the subsequence, $x^{t_j} \to \bar{x}$, meaning that
\begin{equation*}
a^{t_j} = \|x^{t_j}-\bar{x}\| \to 0.
\end{equation*}
Hence $\ell=0$, and therefore
\begin{equation*}
\| x^t-\bar{x} \| \to 0.
\end{equation*}
Thus $x^t \to \bar{x} \in \argmin_{x \in \calC}f(x)$, proving the proposition.

\section{Proof of Lemma~\ref{l:propertiesF}}
\label{sec:prooflempropF}

Property~\ref{i:xones} follows immediately from the fact that $D^\top \ones = 0$, meaning that both $F(x+\alpha \ones)$ and $F(x) + \alpha \ones$ are equal to the output of the circuit, given the input $x + \alpha \ones$. Since this output is unique, we have $F(x + \alpha \ones) = F(x) + \alpha \ones$. 

Next, we prove Property~\ref{i:xy}. Let $x,y \in \R^{N_I}$ be such that $x \preceq y$. Define $z(t) := x+t(y-x)$, where $t \in [0,1]$. Let $p_{HO}(t) \in \R^{|\calV_H| + N_O}$ be the unique vector satisfying
\begin{equation}
\label{zpOH}
\begin{bmatrix}
D_H \\ D_O
\end{bmatrix}\nabla G(D_I^\top z(t) + \begin{bmatrix} D_H^\top & D_O^\top \end{bmatrix} p_{HO}(t)) = 0,
\end{equation}
for all $t \in [0,1]$. It follows by connectedness of $\calG$, the properties of the co-content functions in \ref{i:func}, and the implicit function theorem \cite[Thm. 9.28]{Rudin1976} that $p_{HO}$ is continuously differentiable. 

Due to Property~\ref{i:func}, $\nabla^2 G(D^\top p(t))$ is a diagonal positive definite matrix for all $t \in [0,1]$. Now, let $t \in [0,1]$ and define the Laplacian matrix
\begin{equation*}
L(t) \! := \! D \nabla^2 G(D^\top p(t)) D^\top, \! \text{ where } p(t) := \! \begin{bmatrix}
z(t) \\ p_{HO}(t) 
\end{bmatrix}.
\end{equation*}
We partition the Laplacian as
\begin{equation*}
L(t) = \begin{bmatrix}
L_{11}(t) & L_{12}(t) \\ L_{21}(t) & L_{22}(t) 
\end{bmatrix},
\end{equation*}
where 
\begin{itemize}
\item $L_{11}(t) \in \R^{N_I \times N_I}$,
\item $L_{12}(t) = L_{21}(t)^\top \in \R^{N_I \times (|\calV_H| + N_O)}$, and
\item $L_{22}(t) \in \R^{(|\calV_H| + N_O) \times (|\calV_H| + N_O)}$.
\end{itemize}
In view of \eqref{zpOH}, we obtain
\begin{equation}
\label{Lzp}
\begin{aligned}
0 &= \frac{d}{dt} \begin{bmatrix}
D_H \\ D_O
\end{bmatrix}\nabla G(D^\top p(t)) \\
&= \begin{bmatrix}
D_H \\ D_O
\end{bmatrix}\nabla^2 G(D^\top p(t)) D^\top p'(t) \\
&= L_{21}(t) z'(t) + L_{22}(t) p_{HO}'(t).
\end{aligned}
\end{equation}
Now, Since $\calG$ is connected, the matrix $\begin{bmatrix}
D_H \\ D_O
\end{bmatrix}$ has full row rank, and thus $L_{22}(t)$ is an invertible $M$-matrix. It follows from \cite[Ch. 6, Fact N38]{Berman1994} that $L_{22}^{-1}(t)$ only has nonnegative entries. Furthermore, since $L_{21}(t)$ is an off-diagonal block of the Laplacian matrix, it contains only nonpositive entries. Hence, the matrix $-L_{22}^{-1}(t)L_{21}(t)$ is entry-wise nonnegative. Now, since $z'(t) = y-x$, we obtain by \eqref{Lzp} that
\begin{equation*}
p_{HO}'(t) = - L_{22}^{-1}(t) L_{21}(t) (y-x).
\end{equation*}
Integrating this expression from $t = 0$ to $t=1$ gives
\begin{equation*}
p_{HO}(1) - p_{HO}(0) = -  L_{22}^{-1}(t) L_{21}(t) (y-x).
\end{equation*}
Using entry-wise nonnegativity of both $-L_{22}^{-1}(t)L_{21}(t)$ and $y-x$ yields $p_{HO}(0) \preceq p_{HO}(1)$. In particular, by focusing on the last $N_O$ entries, and using the facts that $z(0) = x$ and $z(1) = y$, we obtain 
\begin{equation*}
F(x) \preceq F(y).
\end{equation*}
This proves the lemma. 

\section{Proof of Lemma~\ref{l:sminmax}}
\label{sec:prooflemsminmax}

We will only show the proof for $\smax_\lambda$. The proof for $\smin_\lambda$ follows analogously. To prove the first statement of the lemma, let $\lambda > 0$ and $x \in \R^n$. Define $M := \max_{i \in [n]} x_i$. Then, obviously, 
\begin{equation*}
\frac{1}{n} \sum_{i = 1}^n e^{\lambda x_i} \leq  e^{\lambda M}.
\end{equation*}
By taking natural logarithms, $\smax_\lambda(x) \leq M$. On the other hand, since $M = x_i$ for some index $i \in [n]$,
\begin{equation*}
\frac{1}{n} \sum_{i = 1}^n e^{\lambda x_i} \geq \frac{1}{n}e^{\lambda M},
\end{equation*}
and thus, 
\begin{equation*}
\smax_\lambda(x) \geq \frac{1}{\lambda} (\lambda M - \ln(n)) = M - \frac{\ln(n)}{\lambda}.
\end{equation*}
We conclude that
\begin{equation*}
M - \frac{\ln(n)}{\lambda} \leq \smax_\lambda(x) \leq M,
\end{equation*}
showing that 
$$
|\smax_\lambda(x) - M | \leq \frac{\ln(n)}{\lambda}. 
$$ 
Now, to prove \eqref{smaxLC}, let $x,y \in \R^n$ and define $d := \|x-y\|_\infty$. Then, obviously, $y-d \ones \preceq x \preceq y+d \ones$. It follows from the definition of $\smax_\lambda$ that $\smax_\lambda(u) \leq \smax_\lambda(v)$ for all $u,v \in \R^n$ such that $u \preceq v$. Therefore,
\begin{equation*}
\smax_\lambda(y-d \ones) \leq \smax_\lambda(x) \leq \smax_\lambda(y + d \ones).
\end{equation*}
Furthermore, for all $\alpha \in \R$,
\begin{equation*}
\smax_\lambda(y+\alpha \ones) = \frac{1}{\lambda} \ln\left(\frac{1}{n} \sum_{i = 1}^n e^{\lambda (y_i + \alpha)} \right),
\end{equation*}
which is equal to 
\begin{equation*}
\frac{1}{\lambda} \ln \left( e^{\lambda \alpha} \frac{1}{n} \sum_{i = 1}^n e^{\lambda y_i} \right) = \smax_\lambda(y) + \alpha.
\end{equation*}
Hence,
\begin{equation*}
\smax_\lambda(y) - d \leq \smax_\lambda(x) \leq \smax_\lambda(y) + d.
\end{equation*}
This shows that 
\begin{equation*}
| \smax_\lambda(x) - \smax_\lambda(y) | \leq d = \|x-y\|_\infty,
\end{equation*}
proving the lemma. 

\vspace{5pt}

\bibliographystyle{plain}
\bibliography{references}

\end{document}